\documentclass[12pt,reqno]{amsart} 
\usepackage{amssymb}
\usepackage{amsmath}
\usepackage{comment}
\usepackage{amsthm}
\usepackage{caption}
\usepackage{mathtools}
\usepackage{fullpage}
\allowdisplaybreaks
\newtheorem{theorem}{Theorem}
\newtheorem{proposition}[theorem]{Proposition}
\newtheorem{lemma}[theorem]{Lemma}
\newcommand{\EndProof}{\hfill $\square$}

\begin{document} 
\numberwithin{equation}{section}

\title{A Splitting into the Double Cover of SL$(3,\mathbb{R})$} 
\author{Edmund Karasiewicz}
\email{ekarasie@ucsc.edu}
\address{Karasiewicz: Department of Mathematics, University of California Santa Cruz, Santa Cruz CA 95064, USA }
%\thanks{This work was supported by the .}
\begin{abstract} We provide a formula for the splitting of a congruence subgroup of SL$(3,\mathbb{R})$ into the double cover of SL$(3,\mathbb{R})$ in terms of Pl\"{u}cker coordinates and prove that the splitting satisfies a twisted multiplicativity. The existence of this splitting and a formula (in terms of a different set of coordinates) was proved by S.D. Miller in an unpublished note; the formula in terms of Pl\"{u}cker coordinates is advantageous to the computation of the Fourier coefficients of an Eisenstein series on the double cover of SL$(3)$ over $\mathbb{Q}$. The computation of these Fourier coefficients will be addressed in a forthcoming work.

\end{abstract}
\maketitle

%%%%%%%%%%%%%%%%%%%%%%%%
%%%%%%INTRO
%%%%%%%%%%%%%%%%%%%%%%%%
\section{Introduction} 

We provide a formula the splitting map (see Subsection \ref{sec:BlockSplit}) of a congruence subgroup of SL$(3,\mathbb{R})$ into the double cover of SL$(3,\mathbb{R})$ in terms of Pl\"{u}cker coordinates. The significance of this map for the number theorist, aside from the arithmetic significance associated with a formula built out of Kronecker symbols, is related to its appearance in the computation of Fourier coefficients of a minimal parabolic Eisenstein series on the double cover of SL$(3)$ over $\mathbb{Q}$.

The existence of this splitting and a formula in terms of a different set of coordinates (see Subsection \ref{sec: BlockParams}) was established by S.D. Miller \cite{M06} in an unpublished note using the work of Banks-Levy-Sepanski \cite{BLS99}; unfortunately, this formula does not appear conducive to the computation of the Fourier coefficients; however, while considering a similar computation over number fields containing the $4$-th roots of unity, Brubaker-Bump-Friedberg-Hoffstein \cite{BBFH07} successfully computed the non-degenerate Fourier coefficients of a minimal parabolic Eisenstein series using a formula for an analog of the splitting in terms of Pl\"{u}cker coordinates. These coordinates are well suited to the computation of the exponential sums and so our main objective is a proof that Miller's formula for the splitting can be written in terms of Pl\"{u}cker coordinates. Additionally, we show that the splitting satisfies a twisted multiplicativity. With this new formula for the splitting we can compute the Fourier coefficients of the Eisenstein series on the double cover of SL$(3)$ over $\mathbb{Q}$. The computation of these Fourier coefficients will be addressed in a forthcoming work. 

Now we will state the main theorem of this paper. The reader can refer to Subsection \ref{sec:PluckSplit} for notation.\newpage

\begin{theorem}\label{IntroSplit} Let $\gamma\in\Gamma_{1}(4)$. Suppose that the Pl\"{u}cker coordinates of $\gamma$ are given by $(4A_1,4B_1,C_1,4A_2,4B_2,C_2)$, such that $A_1>0$, $A_2\neq0$,  and $A_2$ is divisible by fewer powers of $2$ than $A_1$ is. Let $D=(A_1,A_2)$, $D_1 = (D,B_1)$, $D_2=D/D_1$, and let $\epsilon = \left(\frac{-1}{-B_1/D_1}\right)$. Then 
\begin{equation} \label{plucksplitting}
s(\gamma)=\left(\frac{-\epsilon}{-A_1A_2}\right)\left(\frac{A_1/D}{A_2/D}\right)
\left(\frac{B_1/D_1}{A_1/D}\right)\left(\frac{4B_2/D_2}{\operatorname{sign}(A_2)A_2/D}\right)\left(\frac{D_1}{C_1}\right)\left(\frac{D_2}{C_2}\right).\end{equation}
\end{theorem} 

The method of proof is as follows: We use Miller's formula to study certain symmetries satisfied by the splitting map. Using these symmetries we can reduce the computation of the splitting on a general element of the domain to an element with more favorable arithmetic properties. On this special subset of points Miller's formula can be expressed in terms of Pl\"{u}cker coordinates by a direct  (but lengthy) computation.

%In section \ref{chap:Split} we study the splitting $s$. The main result of this section provides a formula for the splitting in terms of Pl\"{u}cker Coordinates as seen in Theorem \ref{IntroSplit}. This result is achieved by studying the effect of certain symmetries on $s$. Additionally, it is established that $s$ satisfies a type of twisted multiplicativity. The formula for the splitting in terms of Pl\"{u}cker coordinates provides us a means to execute the computation of the Fourier coefficients  of the metaplectic Eisenstein series. However, these computations will be addressed in a forthcoming work.

\section{Notation and Conventions} \label{notation}\label{sec:Notation}\label{sec:Comp}
\subsection{SL($3,\mathbb{R}$) and $\widetilde{\text{SL}}(3,\mathbb{R})$}\label{ssec:LieGroups}

This section contains the notation and basic computations that will be used throughout the remainder of this paper. Most of the material contained in this section is fairly standard; Subsection \ref{ssec:2-cocycle} provides one notable exception. This subsection introduces the Banks-Levy-Sepanski 2-cocycle \cite{BLS99} and collects some basic properties of the 2-cocycle. The choice of a 2-cocycle defines a multiplication on $\widetilde{\text{SL}}(3,\mathbb{R})$ and thus plays an important role in the determination of the splitting $s$.

Let $\widetilde{G}=\widetilde{\text{SL}}(3,\mathbb{R})$, the double cover of $G=\text{SL}(3,\mathbb{R})$. Let $N_{R}$ be equal to the group of $3\times3$ upper triangular unipotent matrices, let $B_{R}$ be the group of $3\times3$ upper triangular matrices with, and let $T_R$ be the subgroup of diagonal matrices. In each case the subgroups entries in the ring $R$; when $R=\mathbb{R}$ we may omit the subscript. Let $A$ be the subgroup of $T$ with positive entries. Finally we will use the following discrete subgroups:\begin{equation}\label{gamma14+infty}\begin{array}{lll}
\Gamma=\Gamma_{1}(4)= 
\{\gamma\in \text{SL}(3,\mathbb{Z})|\gamma \equiv \left(\begin{smallmatrix}
1&*&*\\
0&1&*\\
0&0&1\end{smallmatrix}
\right)\text{(mod }4)\} &, & \Gamma_{\infty}=N_{\mathbb{R}}\cap\text{SL}(3,\mathbb{Z}).
\end{array}\end{equation}
We will use the following notation for elements of $T$ and $N$ respectively: $
\text{t}(a,b,c) = \left(\begin{smallmatrix}
a & & \\
 & b & \\
  & & c
\end{smallmatrix}\right)$, and 
$n(x,y,z) = \left(\begin{smallmatrix}
1 & x & z\\
0 & 1 & y\\
0 & 0 & 1\end{smallmatrix}
\right)
$.

As a set $\widetilde{\text{SL}}(3,\mathbb{R})\cong \text{SL}(3,\mathbb{R})\times\{\pm1\}$. As a topological space $\widetilde{\text{SL}}(3,\mathbb{R})$ is given the nontrivial covering space topology whose existence is guaranteed by the fact that the fundamental group of SL$(3,\mathbb{R})$ is isomorphic to $\mathbb{Z}/2\mathbb{Z}$.

The Banks-Levy-Sepanski 2-cocycle $\sigma:\text{SL}(3,\mathbb{R})\times\text{SL}(3,\mathbb{R})\rightarrow \{\pm1\}$, constructed in \cite{BLS99} and recalled in Subsection \ref{ssec:2-cocycle}, defines the group multiplication as follows:
\begin{equation}\label{metamult}(g_1,\epsilon_1)(g_2,\epsilon_2) = (g_1g_2,\epsilon_1\epsilon_2\sigma(g_1,g_2)).
\end{equation}

The following list establishes some notation for elements of the Weyl group of SL$(3,\mathbb{R})$:
\[\begin{array}{llllll}
w_{\alpha_1}&=&\left(\begin{smallmatrix}
0 & -1 & 0\\
1 & 0 & 0\\
0 & 0 & 1\end{smallmatrix}\right),&
w_{\alpha_2}&=&\left(\begin{smallmatrix}
1 & 0 & 0\\
0 & 0 & -1\\
0 & 1 & 0\end{smallmatrix}\right),\\
w_{\alpha_1}w_{\alpha_2}&=&\left(\begin{smallmatrix}
0 & 0 & 1\\
1 & 0 & 0\\
0 & 1 & 0\end{smallmatrix}\right),&
w_{\alpha_2}w_{\alpha_1}&=&\left(\begin{smallmatrix}
0 & -1 & 0\\
0 & 0 & -1\\
1 & 0 & 0\end{smallmatrix}\right),
\end{array}\]
\[\text{and} \begin{array}{lllllllll}
w_{\ell}&=&w_{\alpha_1}w_{\alpha_2}w_{\alpha_1}&=&w_{\alpha_2}w_{\alpha_1}w_{\alpha_2}&=&\left(\begin{smallmatrix}
0 & 0 & 1\\
0 & -1 & 0\\
1 & 0 & 0\end{smallmatrix}\right).& &
\end{array}\]
The representatives of the Weyl group listed above are those defined in Section 3 of \cite{BLS99}. These representatives are used in the formula for the 2-cocycle, which can be found in Section 4 of the previously cited paper. 

In \cite{M06}, Miller constructs a map $S:\Gamma_{1}(4)\hookrightarrow\widetilde{\text{SL}}(3,\mathbb{R})$ such that $S(\gamma) = (\gamma,s(\gamma))$, where $s(\gamma)\in \{\pm1\}$. This map $S$ is a splitting of $\Gamma_{1}(4)$ into $\widetilde{SL}(3,\mathbb{R})$. The map $s$, which by abuse of terminology will also be called the splitting, is defined in section \ref{sec:BlockSplit}. 

\subsection{Pl\"{u}cker Coordinates} \label{sec:SL3R}

In this section let $\Gamma =$ SL$(3,\mathbb{Z})$. The \textbf{Pl\"{u}cker coordinates} of $\Gamma_{\infty}\backslash \Gamma$ are defined presently. Given $\left( \begin{smallmatrix}
a & b & c \\
d & e & f \\
g & h & i \end{smallmatrix} \right)\in $ SL$(3,\mathbb{R})$, define six parameters as follows:
\begin{equation}\label{PluckerCoordsAgain}\begin{array}{ccc}
A_1^\prime = -g, & B_1^\prime = -h, & C_1^\prime = -i\\
A_2^\prime = -(dh-eg), & B_2^\prime = (di-fg), & C_2^\prime = -(ei-fh).\end{array}
\end{equation}
\begin{theorem} \label{theorem:PluckerCoords}(\cite[Ch 5]{BGL3}) The map taking $\left( \begin{smallmatrix}
a & b & c \\
d & e & f \\
g & h & i \end{smallmatrix} \right)\in $ SL$(3,\mathbb{R})$ to $(A_1^\prime,B_1^\prime,C_1^\prime,A_2^\prime,B_2^\prime,C_2^\prime)$ defines a bijection between the coset space $N\backslash$SL$(3,\mathbb{\mathbb{R}})$ and the set of all \newline$(A_1^\prime,B_1^\prime,C_1^\prime,A_2^\prime,B_2^\prime,C_2^\prime)\in \mathbb{R}^6$ such that: $A_1^\prime C_2^\prime+B_1^\prime B_2^\prime+C_1^\prime A_2^\prime= 0$, not all of $A_1^\prime,B_1^\prime,C_1^\prime$ equal $0$, and not all of $A_2^\prime,B_2^\prime,C_2^\prime$ equal $0$.  Furthermore, a coset in $N\backslash$SL$(3,\mathbb{\mathbb{R}})$ contains an element of $\Gamma$ if and only if $A_1^\prime,B_1^\prime,C_1^\prime$ are coprime integers and $A_2^\prime,B_2^\prime,C_2^\prime$ are coprime integers.\end{theorem}

Versions of this result hold for other congruence subgroups. Let
$A_1^\prime = 4A_1$,  $A_2^\prime = 4A_2$,
$B_1^\prime = 4B_1$,  $B_2^\prime =4B_2$,
$C_1^\prime = C_1$,  $C_2^\prime = C_2$.
The coset space 
 $\Gamma_{\infty}\backslash \Gamma_1(4)$ can be identified with
\begin{equation}\label{cosetparam2}\left\{\begin{array}{c}
(4A_1,4B_1,C_1,4A_2,4B_2,C_2)\in \mathbb{Z}^6 | A_1C_2+4B_1B_2+C_1A_2 = 0,\\
(A_i,B_i,C_i) = 1,\,C_j\equiv -1\text{ (mod } 4)
\end{array}\right\}
\begin{array}{l}
\text{}\\
.\end{array}
\end{equation}

%%%%%%%%%%%%%%%%
%%%%%COMPUTATIONS
%%%%%%%%%%%%%%%%

\subsection{Kronecker Symbol} \label{ssec:Kronecker}

If $a,b\in\mathbb{R}^{\times}$, then the Hilbert Symbol, $(a,b)_{\mathbb{R}}$, is to be $1$ if $a$ or $b$ is positive and $-1$ if $a$ and $b$ are negative. 
Let $n\in \mathbb{Z}_{\neq 0}$ with prime factorization $n=\epsilon p_{1}^{e_1}\ldots p_{\ell}^{e_\ell}$, where $\epsilon=\pm1$. If $k\in\mathbb{Z}$ the Kronecker Symbol is defined by $\big(\frac{k}{n}\big) = \big(\frac{k}{\epsilon}\big)\big(\frac{k}{p_{1}}\big)^{e_1}\ldots\big(\frac{k}{p_{\ell}}\big)^{e_{\ell}}$, where $\big(\frac{k}{p_{i}}\big)$ is the Legendre symbol when $p_i$ is an odd prime, $\big(\frac{k}{\epsilon}\big)=(k,\epsilon)_{\mathbb{R}}$, and $\big(\frac{k}{2}\big)$ is equal to $0$ if $2$ divides $k$, $1$ if $k\equiv \pm1 (\text{mod } 8)$, and $-1$ if $k\equiv \pm3 (\text{mod } 8)$. The formula can be extended to $n=0$ by setting $\big(\frac{k}{0}\big)$ equal to $1$ if $k=\pm1$, and $0$ otherwise.

The following proposition collects some facts about the Kronecker Symbol. These results can be found in \cite{IK04}.
\begin{proposition} \label{Kronecker}  (Properties of the Kronecker Symbol) Let $a,b,m,n\in\mathbb{Z}$, $\epsilon=\pm1$, and let $n^\prime$ and $m^\prime$ denote the odd part of $n$ and $m$, respectively.
\begin{enumerate}
\item If $ab\neq0$, then $\left(\frac{a}{n}\right)\left(\frac{b}{n}\right)=\left(\frac{ab}{n}\right).$
\item If $mn\neq0$, then $\left(\frac{a}{m}\right)\left(\frac{a}{n}\right)=\left(\frac{a}{mn}\right).$
\item Let $m$ be equal to $4n$ if $n\equiv 2\, (\text{mod }4)$, otherwise let $m$ be equal to $n$. If $n>0$ and $a\equiv b\,(\text{mod }m)$, then $\left(\frac{a}{n}\right)=\left(\frac{b}{n}\right)$.
\item Let $b$ be equal to $4|a|$ if $a\equiv 2\, (\text{mod }4)$, otherwise let $b$ be equal to $|a|$. If $a\not\equiv 3\,(\text{mod }4)$ and $m\equiv n\,(\text{mod } b)$, then $(\frac{a}{m})=(\frac{a}{n})$.
\item $\left(\frac{-1}{n}\right)=(-1)^{\frac{n^\prime-1}{2}}$ and $\left(\frac{2}{n^\prime}\right)=(-1)^{\frac{(n^\prime)^2-1}{8}}$.
\item (Quadratic Reciprocity) If gcd$(m,n) =1$, then $\big(\frac{m}{n}\big)\big(\frac{n}{m}\big)=(n,m)_{\mathbb{R}}(-1)^{\frac{(m^\prime-1)(n^\prime-1)}{4}}$.
\item$(\frac{(\frac{-1}{m})}{n}) = (-1)^{(\frac{m^\prime-1}{2})(\frac{n^\prime-1}{2})}=((\frac{-1}{m}),(\frac{-1}{n}))_{\mathbb{R}}$.
\end{enumerate}
\end{proposition}

\subsection{2-cocycle}\label{ssec:2-cocycle}

This section begins with a formula for the Banks-Levy-Sepanski 2-cocycle \cite{BLS99} as presented in Miller \cite{M06} and then goes on to collect some facts about the 2-cocycle that will be used in subsequent computations.

Let $g\in $ SL$(3,\mathbb{R})$ with Pl\"{u}cker coordinates $(A_1^\prime,B_1^\prime,C_1^\prime,A_2^\prime,B_2^\prime,C_2^\prime)$. Let $X_1(g)=\text{det}(g)$, let $X_2(g)$ be the first non-zero element of the list $-A_2,B_2,-C_2$, let $X_3(g)$ be the first nonzero element of the the list $-A_1,-B_1,-C_1$, and let 
$\Delta(g)=\bigg(\begin{smallmatrix}
X_1(g)/X_2(g)& 0& 0\\
0 & X_2(g)/X_3(g) & 0\\
0 & 0 & X_3(g)
\end{smallmatrix}\bigg)$.

If $g_1,g_2\in G$ such that $g_1=naw_{1}\ldots w_{k}n^\prime$ is the Bruhat decomposition of $g_1$, then in Section 4 of \cite{BLS99}, Banks-Levy-Sepanski show that the 2-cocycle $\sigma$ satisfies the formula
\begin{equation}\label{BLSdef}\sigma(g_1,g_2) = \sigma(a,w_{1}\ldots w_{k}n^\prime g_2)\sigma(w_{1},w_{2}\ldots w_{k}n^\prime g_2)\ldots\sigma(w_{k-1},w_{k}n^\prime g_2)\sigma(w_{k},n^\prime g_2),\end{equation}
where each factor can be computed using the following rules: Let $h\in G$, $a\in T$, then
\begin{align*}
\sigma(t(a_1,a_2,a_3),t(b_1,b_2,b_3))=& (a_1,b_2)(a_1,b_3)(a_2,b_3),\\
\sigma(a,h) = \sigma(a,\Delta(h)),\text{ and }
\sigma(w_\alpha,h) =& \sigma(\Delta(w_\alpha h)\Delta(h),-\Delta(h)),
\end{align*} 

The next lemma collects some simple identities involving $\sigma$.
\begin{lemma} \label{lemma:cocycle}(Banks-Levy-Sepanski \cite{BLS99})
Let $n,n_1,n_2\in N$ and let $a\in A$, Then\\
$\sigma(n_1g_1,g_2n_2)=\sigma(g_1,g_2)$,
$\sigma(g_1n,g_2) = \sigma(g_1,ng_2)$,
$\sigma(n,g) = \sigma(g,n) = 1$, and
$\sigma(g,a)=1$.

\end{lemma}

\textbf{Proof:} The first three identities follow from Lemma 4 in Section 3 of \cite{BLS99}. The fourth identity follows as the Steinberg symbol may be taken to be the Hilbert symbol, and the combination of Theorem 7a) in \cite{BLS99} and the second to last equation in Section 4 of \cite{BLS99}.\EndProof\\

%%%%%%%%%%%%%%%%%%%%%%%
%%%%%%%SPLITTING
%%%%%%%%%%%%%%%%%%%%%%%

\section{The Splitting} \label{chap:Split}

This section contains a study of the splitting $\Gamma_{1}(4)\rightarrow\widetilde{\text{SL}}(3,\mathbb{R}),$ where $\gamma\mapsto(\gamma,s(\gamma))$ and $s:\Gamma_1(4)\rightarrow\{\pm1\}$ is the map defined on line (\ref{blocksplitting}) in Subsection \ref{sec:BlockSplit}, culminating in the proof of our main result, Theorem \ref{IntroSplit}. 

%%%%%%%%%%%%%%%%%%%%%%%%%%%%%%%%%%
%%%%%%%%%Block Parameters
%%%%%%%%%%%%%%%%%%%%%%%%%%%%%%%%%%

\subsection{Block Parameters} \label{sec: BlockParams}
Any $\gamma\in \Gamma$ can be written as 
\begin{equation} \label{blockfactorization}
\gamma = n
\left( \begin{smallmatrix}
a_1 & b_1 & 0\\
c_1 & d_1 & 0\\
0 & 0 & 1 \end{smallmatrix} \right)
\left( \begin{smallmatrix}
a_2 & 0 & b_2\\
0 & 1 & 0\\
c_2 & 0 & d_2 \end{smallmatrix} \right)\
\left( \begin{smallmatrix}
1 & 0 & 0\\
0 & a_3 & b_3\\
0 & c_3 & d_3 \end{smallmatrix} \right).
\end{equation}  
where $n\in N$ and all of the other factors are in $\Gamma$. Note that this factorization is non-unique. The parameters $a_i,b_i,c_i,d_i$ will be called \textbf{block parameters}. In \cite{M06}, Miller proves that $\Gamma_{1}(4)$ is maximal in SL$(3,\mathbb{Z})$ such that there exists a splitting homomorphism $\Gamma_{1}(4)\rightarrow\widetilde{\text{SL}}(3,\mathbb{R})$ and provides a formula for this splitting homomorphism using equation (\ref{blockfactorization}) and the 2-cocycle of Banks-Levy-Sepanski \cite{BLS99}.

Theorem \ref{theorem:PluckerCoords} (above) shows how Pl\"{u}cker coordinates can be used to provide a description of $\Gamma_{\infty}\backslash\Gamma$. By relating the block parameters and the Pl\"{u}cker coordinates, this description can be used to identify a unique matrix representative of $\Gamma_{\infty}\backslash\Gamma$ in terms of the block parameters. Proposition \ref{blockPluck} below makes this relationship between block parameters and Pl\"{u}cker coordinates explicit, and Proposition \ref{blockrep} shows how this relationship can be used to pick a unique representative of $\Gamma_{\infty}\backslash\Gamma$ using the block parameters.

%Prove factorization.  This factorization comes from Miller's GftM.
\begin{proposition} \label{blockPluck}If $\gamma = n
\left( \begin{smallmatrix}
a_1 & b_1 & 0\\
c_1 & d_1 & 0\\
0 & 0 & 1 \end{smallmatrix} \right)
\left( \begin{smallmatrix}
a_2 & 0 & b_2\\
0 & 1 & 0\\
c_2 & 0 & d_2 \end{smallmatrix} \right)
\left( \begin{smallmatrix}
1 & 0 & 0\\
0 & a_3 & b_3\\
0 & c_3 & d_3 \end{smallmatrix} \right)\in\Gamma_{1}(4)$,
then $\gamma$ has Pl\"{u}cker coordinates given by 
\begin{equation}\label{blocktopluck}
\begin{array}{lll}A_1^\prime = -c_2,& B_1^\prime = -d_2c_3,& C_1^\prime=-d_2d_3\\
A_2^\prime = -(c_1c_3-d_1c_2a_3),&B_2^\prime = c_1d_3-d_1c_2b_3,&C_2^\prime = -d_1d_2.\end{array}\end{equation}
\end{proposition} 

\noindent\textbf{Proof:} This follows directly from the definitions.\EndProof
%NOTE: As Miller's parameters are coordinates of matrices the equation from Bump's parameters is satisfied.  

\begin{proposition}  \label{blockrep}The coset of $\Gamma_{\infty}\backslash\Gamma$ with Pl\"{u}cker coordinates $(A_1^\prime,B_1^\prime,C_1^\prime,A_2^\prime,B_2^\prime,C_2^\prime)$ is represented by the matrix 
\[\left( \begin{smallmatrix}
a_1 & b_1 & 0\\
c_1 & \frac{-C_2^\prime}{(B_1^\prime,C_1^\prime)\epsilon} & 0\\
0 & 0 & 1 \end{smallmatrix} \right)
\left( \begin{smallmatrix}
a_2 & 0 & b_2\\
0 & 1 & 0\\
-A_1^\prime & 0 & (B_1^\prime,C_1^\prime)\epsilon \end{smallmatrix} \right)\
\left( \begin{smallmatrix}
1 & 0 & 0\\
0 & a_3 & b_3\\
0 & \frac{-B_1^\prime}{(B_1^\prime,C_1^\prime)\epsilon} & \frac{-C_1^\prime}{(B_1^\prime,C_1^\prime)\epsilon} \end{smallmatrix} \right),\]
where \begin{enumerate}
\item $\epsilon = \left(\frac{-1}{(B_1^\prime,C_1^\prime)}\right).$
\item $a_3$ is the smallest positive integer satisfying $a_3\equiv d_3^{-1} \text{(mod }c_3)$, such that $\frac{d_1c_2a_3-A_2^\prime}{c_3}$ $\in4\mathbb{Z}$.
\item $a_i$ ($i\neq 3$) is the smallest positive integer satisfying $a_i\equiv d_i^{-1} \text{(mod }c_i)$. 
\item $c_1=\frac{d_1c_2a_3-A_2^\prime}{c_3}$.
\item $b_j=\frac{a_jd_j-1}{c_j}$ for $j=1,2,3$.
\end{enumerate}
\end{proposition}
The appearance of `smallest' in items two and three is a choice made to ensure uniqueness.

\noindent\textbf{Proof:} Let $\gamma\in \Gamma$ with Pl\"{u}cker coordinates $(A_1^\prime,B_1^\prime,C_1^\prime,A_2^\prime,B_2^\prime,C_2^\prime)$, and say 
\[\gamma = n
\left( \begin{smallmatrix}
a_1 & b_1 & 0\\
c_1 & d_1 & 0\\
0 & 0 & 1 \end{smallmatrix} \right)
\left( \begin{smallmatrix}
a_2 & 0 & b_2\\
0 & 1 & 0\\
c_2 & 0 & d_2 \end{smallmatrix} \right)\
\left( \begin{smallmatrix}
1 & 0 & 0\\
0 & a_3 & b_3\\
0 & c_3 & d_3 \end{smallmatrix} \right).\]  
Using line (\ref{blocktopluck}), the block parameters may be expressed as follows:
\[c_2 = -A_1^\prime,\hspace{5mm} d_2=\left(\frac{-1}{(B_1^\prime,C_1^\prime)}\right)(B_1^\prime,C_1^\prime), \hspace{5mm}c_3 = \frac{-B_1^\prime}{d_2}, \hspace{5mm}d_3 = \frac{-C_1^\prime}{d_2}\]
\[d_1 = \frac{-C_2^\prime}{d_2}, \hspace{5mm} c_1\frac{-B_1^\prime}{d_2}\equiv -A_2^\prime\, (\text{mod} \frac{A_1^\prime C_2^\prime}{d_2} = d_1c_2).\]

This determines $d_1, c_2, d_2, c_3, d_3\in \mathbb{Z}$ such that $(c_2,d_2) = (c_3,d_3) = 1$, $c_2, c_3\equiv 0 \text{(mod } 4)$, and $d_1, d_2, d_3 \equiv 1\text{(mod } 4)$. Consider the expression for $A_2^\prime$ in line (\ref{blocktopluck}). As $(c_3,d_1c_2)$ divides $A_2^\prime$, there are infinitely many pairs of integers $c_1$ and $a_3$ that satisfy that equation. Take $a_3$ to be the smallest positive integer satisfying $a_3\equiv d_3^{-1} \text{(mod }c_3)$ and $c_1\in4\mathbb{Z}$. The relation $(A_2^\prime,B_2^\prime,C_2^\prime)=1$ implies that $(c_1,d_1)=1$. Similarly, take $a_i$ ($i\neq 3$) to be the smallest positive integer satisfying $a_i\equiv d_i^{-1} \text{(mod }c_i)$.
%NOTE: We get a_3=1 mod 4 for free as 4|c_3 and d_3=1 mod 4.

Define $b_i$ to be $\frac{a_id_i-1}{c_i}$. Finally we can show that this choice for $b_3$ is consistent with the expression for $B_2^\prime$.\EndProof

%NOTE:  A linear Diophantine equation Ax+By = C has a solution if and only if (A,B) divides C.  If (x,y) is such a solution, then all possible solutions are of the form (x+nB/(A,B),y-nA/(A,B)).  Thus we may find a unique solution by stipulating that x\in \{0,1,\ldots,B/(A,B)-1\}.

%%%%%%%%%%%%%%%%%%%%%%%%%%%%%%%%%%%%
%%%%%%%Symmetries of Bump's Parameters
%%%%%%%%%%%%%%%%%%%%%%%%%%%%%%%%%%%%

\subsection{Symmetries of Pl\"{u}cker Coordinates} \label{sec:PluckerSym}

The current section collects the symmetries of the splitting map required to prove Theorem \ref{IntroSplit}.
\begin{proposition} \label{prop:PluckerSym}
Let $M\in \text{SL}(3,\mathbb{R})$ with Pl\"{u}cker coordinates \newline$(4A_1,4B_1,C_1,4A_2,4B_2,C_2)$, $n=n(x,y,z)\in N$, $S_2=\text{t}(1,-1,1)$, and $S_3 = \text{t}(1,1,-1)$. Then:
\begin{enumerate}
\item The matrix $nMn^{-1}$ has Pl\"{u}cker coordinates \[(4A_1,4B_1-4A_1x,C_1-4B_1y+4A_1(xy-z),4A_2,4B_2+4A_2y,C_2+4B_2x+4A_2z).\]
\item The matrix $S_3MS_3$ has Pl\"{u}cker coordinates $(-4A_1,-4B_1,C_1,-4A_2,4B_2,C_2)$.
\item The matrix $S_2MS_2$ has Pl\"{u}cker coordinates $(4A_1,-4B_1,C_1,4A_2,-4B_2,C_2)$.
\item The matrix $w_{\ell} M^{-t}w_{\ell}$ has Pl\"{u}cker coordinates $(4A_2,-4B_2,C_2,4A_1,-4B_1,C_1)$.
\item Let $M\in\Gamma_{1}(4)$. If $D$ divides $(A_1,A_2)$, $D_1=(D,B_1)$, $D=D_1D_2$, and $T=\text{t}(1,D_2^{-1},D^{-1})$ then $TMT^{-1}\in\text{SL}(3,\mathbb{Z})$ has Pl\"{u}cker coordinates 
\[(4A_1/D,4B_1/D_1,C_1,4A_2/D,(4B_2)/D_2,C_2).\]
Furthermore, $TMT^{-1}\in\Gamma_1(4)$ if and only if $D_2$ divides $B_2$. 
\end{enumerate}
\end{proposition}
The proof is straightforward matrix algebra and will be omitted. The last symmetry is the most important in terms of the reduction step. Many of the technical aspects of the proof of Theorem \ref{IntroSplit} result from the failure of $D_2$ to divide $B_2$ in general. 

%%%%%%%%%%%%%%%%%%%%%%%%%%%%%%%%%%%%
%%%%%%%Some Results on Double Cosets
%%%%%%%%%%%%%%%%%%%%%%%%%%%%%%%%%%%%

\subsection{Some Results on Double Coset Spaces} \label{sec:DoubleCoset}

This section describes some structural features of the set $\Gamma_\infty \backslash\Gamma_1(4)/\Gamma_\infty$.  This double coset space is important for two reasons. First, the map $s:\Gamma\rightarrow \{\pm1\}$ is left and right $\Gamma_{\infty}$-invariant and so it naturally descends to a map with domain $\Gamma_\infty \backslash\Gamma_1(4)/\Gamma_\infty$. Second, certain subsets of $\Gamma_\infty \backslash\Gamma_1(4)/\Gamma_\infty$ index exponential sums appearing in the computation of the Fourier coefficients of the metaplectic Eisenstein series. This section includes an explicit description of the sets indexing the aforementioned exponential sums and proves that the sets indexing the exponential sums exhibit a multiplicative structure. 

Let 
\begin{equation} \label{BoldSA1A2}
\mathbb{S}(A_1,A_2)=\{\gamma\in\Gamma_1(4)| \gamma \text{ has Pl\"{u}cker coordinates of the form }(4A_1,*,*,4A_2,*,*)\}.\end{equation} 
These sets are left $\Gamma_\infty$-invariant and $\Gamma_1(4) = \coprod_{(A_1,A_2)\in \mathbb{Z}^2} \mathbb{S}(A_1,A_2).$  Thus, the Pl{\"u}cker coordinates can be passed to the space $\Gamma_\infty\backslash \mathbb{S}(A_1,A_2)$, and $\Gamma_\infty \backslash\Gamma_1(4) = \coprod \Gamma_{\infty}\backslash\mathbb{S}(A_1,A_2).$ The notation and calculations from section \ref{sec:PluckerSym} provide maps $\text{Ad}(S):\mathbb{S}(A_1,A_2)\rightarrow \mathbb{S}(-A_1,-A_2)$ given by $M\mapsto S_{i}MS_{i}^{-1}$, $\phi:\mathbb{S}(A_1,A_2)\rightarrow \mathbb{S}(A_2,A_1)$ given by $M\mapsto w_{\ell} M^{-t}w_{\ell}$, and $\text{Ad}(n):\mathbb{S}(A_1,A_2)\rightarrow \mathbb{S}(A_1,A_2)$ given by $M\mapsto nMn^{-1}$, where $n\in \Gamma_{\infty}$.  The maps descend to maps on the double coset spaces.
%%%%%%%%%%%%%%%%%%%%
%%%%%%%%%%%%%%%%%%%%12/13/16 I have changed the sign of the representatives. Originally there were minus signs. However, I cannot remember why I made that choice the first time around. Perhaps it was due to the exp sums. The new formulas do not include minus signs so I decided to change this definition.
%%%%%%%%%%%%%%%%%%%%
%%%%%%%%%%%%%%%%%%%%

When $A_1,A_2\neq 0$, the next proposition establishes that unique representatives of the double coset space $\Gamma_\infty\backslash \mathbb{S}(A_1,A_2)/\Gamma_\infty$ are given by the elements of
\begin{multline}\label{SA1A2}S(A_1,A_2) \stackrel{\text{def}}{=} \{(4A_1,4B_1,C_1,4A_2,4B_2,C_2)\in \mathbb{Z}^6 |
A_1C_2+4B_1B_2+C_1A_2 = 0,\\
(A_i,B_i,C_i) = 1,\text{ } C_j\equiv -1\text{(mod } 4), \text{ } \frac{B_1}{A_1},\frac{B_2}{A_2},\frac{C_2}{4A_2}\in[0,1)\}.
\end{multline}

\begin{proposition} \label{DoubleCosetProp} Let $A_1,A_2\neq 0$ and let $\gamma\in \mathbb{S}(A_1,A_2)$ with Pl\"{u}cker coordinates \newline$(4A_1,4B_1,C_1,4A_2,4B_2,C_2)$.  Then there is a unique $n\in \Gamma_{\infty}$ such that the Pl\"{u}cker coordinates of $\gamma n$ live in $S(A_1,A_2)$. This induces a bijective map \[\Gamma_\infty\backslash \mathbb{S}(A_1,A_2)/\Gamma_\infty \rightarrow S(A_1,A_2).\]
 
\end{proposition}

\noindent\textbf{Proof:}  The result follows from the description of the action of Ad($n$) on the coordinates in Subsection \ref{sec:PluckerSym}.\EndProof

Now we turn to the main result of this section; the sets $\Gamma_\infty\backslash \mathbb{S}(A_1,A_2)/\Gamma_\infty$ exhibit a multiplicative structure. This result is essentially an application of the Chinese Remainder Theorem. 

\begin{proposition} \label{mult} Let $A_1,\alpha_1> 0$, $A_2,\alpha_2\neq0,$ suppose that $(A_1A_2,\alpha_1\alpha_2) = 1$, $A_1,A_2$ are odd, and suppose that $A_1\alpha_1+A_2\alpha_2\equiv 0$ (mod 4). Let $\mu=(\frac{-1}{-A_1A_2})$. Then 
\[\Gamma_\infty\backslash \mathbb{S}(A_1\alpha_1,A_2\alpha_2)/\Gamma_\infty \cong \Gamma_\infty\backslash \mathbb{S}(A_1,\mu A_2)/\Gamma_\infty \times\Gamma_\infty\backslash \mathbb{S}(\alpha_1,-\mu\alpha_2)/\Gamma_\infty.\]

The bijection is induced by the map
\[(4A_1\alpha_1,4B_1,C_1,4A_2\alpha_2,4B_2,C_2)\mapsto\]
\begin{multline}((4A_1,4B_1,\mathcal{C}_1,\mu4 A_2,4B_2,\gamma C_2),\\(4\alpha_1,4B_1,\Big(\frac{-1}{A_2}\Big) A_2C_1,-\mu4\alpha_2,-\Big(\frac{-1}{A_2}\Big)\mu4B_2,-\mu \Big(\frac{-1}{A_2}\Big)A_1C_2)),
\end{multline}
where:
\begin{enumerate}
\item$\mathcal{C}_1=\frac{-A_1\gamma C_2-4B_1B_2}{\mu A_2}$.
\item$\gamma$ is the smallest positive integer such that $\gamma\equiv 1$ (mod 4) and $\gamma\equiv \alpha_1$ (mod $A_2$).
\end{enumerate}
\end{proposition}
\noindent\textbf{Proof:}  This proof involves the verification of many simple claims. For the sake of space we will state the claims that need to be verified and content ourselves with proving those which are less straightforward. 

To establish this bijection, we will construct a map 
\[\phi:\Gamma_\infty\backslash \mathbb{S}(A_1\alpha_1,A_2\alpha_2)/\Gamma_\infty \rightarrow \Gamma_\infty\backslash \mathbb{S}(A_1,\mu A_2)/\Gamma_\infty \times\Gamma_\infty\backslash \mathbb{S}(\alpha_1,-\mu\alpha_2)/\Gamma_\infty\] 
and a map
\[\psi : \Gamma_\infty\backslash \mathbb{S}(A_1,\mu A_2)/\Gamma_\infty \times\Gamma_\infty\backslash \mathbb{S}(\alpha_1,-\mu\alpha_2)/\Gamma_\infty\rightarrow \Gamma_\infty\backslash \mathbb{S}(A_1\alpha_1,A_2\alpha_2)/\Gamma_\infty\] 
 such that $\phi$ and $\psi$ are inverses. To describe this map and its inverse, we will need an integer $\gamma\in\mathbb{Z}$ such that $\gamma\equiv 1$ (mod 4) and $\gamma\equiv \alpha_1$ (mod $A_2$) ($A_2$ must be odd for this step). To avoid ambiguities we will stipulate that $\gamma$ is the least positive integer satisfying these conditions. It will be useful to write $\gamma=\alpha_1+\ell A_2$ for some integer $\ell$.
 
 We begin with the construction of $\phi$. Let \newline$(4A_1\alpha_1,4B_1,C_1,4A_2\alpha_2,4B_2,C_2)\in S(A_1\alpha_1,A_2\alpha_2)$.  Consider 
\begin{multline}\label{6tuples}((4A_1,4B_1,\mathcal{C}_1,\mu4A_2,4B_2,\gamma C_2),\\(4\alpha_1,4B_1,\left(\frac{-1}{A_2}\right) A_2C_1,-\mu4\alpha_2,-\left(\frac{-1}{A_2}\right)\mu4B_2,-\mu \left(\frac{-1}{A_2}\right)A_1C_2)),\end{multline}
where $\mathcal{C}_1=\frac{-A_1\gamma C_2-4B_1B_2}{\mu A_2}$.

This will be the value of $\phi$ once we adjust for cosets. We must show that this pair of 6-tuples resides in $\Gamma_\infty\backslash \mathbb{S}(A_1,\mu A_2) \times\Gamma_\infty\backslash \mathbb{S}(\alpha_1,-\mu\alpha_2)$.
To show that \newline$(4A_1,4B_1,\mathcal{C}_1,\mu 4A_2,4B_2,\gamma C_2)\in \Gamma_\infty\backslash \mathbb{S}(A_1,\mu A_2)$, we  verify the following six conditions:
\begin{comment}
\begin{enumerate}
\item $A_1\gamma C_2+4B_1B_2+\mathcal{C}_1\mu A_2=0$
\item $\mathcal{C}_1\in\mathbb{Z}$
\item $\mathcal{C}_1\equiv -1$ (mod $4$)
\item $\gamma C_2\equiv-1$ (mod $4$)
\item $(A_1,B_1,\mathcal{C}_1)=1$
\item $(\mu A_2,B_2,\gamma C_2)=1$
\end{enumerate}
\end{comment}
\[\begin{array}{ll}
(1)\hspace{.5cm} A_1\gamma C_2+4B_1B_2+\mathcal{C}_1\mu A_2=0, & (2)\hspace{.5cm} \mathcal{C}_1\in\mathbb{Z}\\
(3) \hspace{.5cm}\mathcal{C}_1\equiv -1 \text{ (mod }4), & (4) \hspace{.5cm} \gamma C_2\equiv-1\text{ (mod }4)\\
(5) \hspace{.5cm} (A_1,B_1,\mathcal{C}_1)=1, & (6) \hspace{.5cm} (\mu A_2,B_2,\gamma C_2)=1.
\end{array}\]

We will prove item (2) from this list. We have
\begin{equation}\begin{array}{rlcc}\mathcal{C}_1=&\frac{-A_1\gamma C_2-4B_1B_2}{\mu A_2}&\hspace{.5cm}&(\text{def of }\mathcal{C}_1)\\
=&\mu \alpha_2\frac{-A_1\alpha_1C_2-4B_1B_2}{A_2\alpha_2}-\mu \ell A_1C_2&\hspace{.5cm}&(\gamma=\alpha_1+\ell A_2)\\
=&\mu\alpha_2C_1-\mu \ell A_1C_2.&\hspace{.5cm}&(A_1\alpha_1C_2+4B_1B_2+A_2\alpha_2C_1=0)
\end{array}
\end{equation}
Therefore $\mathcal{C}_1$ is an integer. For later reference we record the identity
\begin{equation}
\mathcal{C}_1=\mu\alpha_2C_1-\mu \ell A_1C_2.\label{C1identity}
\end{equation}
  
Next to show that \newline$(4\alpha_1,4B_1,\left(\frac{-1}{A_2}\right) A_2C_1,-\mu4\alpha_2,-\left(\frac{-1}{A_2}\right)\mu 4B_2,-\mu \left(\frac{-1}{A_2}\right)A_1C_2))\in\Gamma_\infty\backslash \mathbb{S}(\alpha_1,-\mu\alpha_2)$, we must verify the following five conditions:
\begin{comment}
\begin{enumerate}
\item  $\alpha_1(-\mu\left(\frac{-1}{A_2}\right)A_1C_2)+4(B_1)(-\mu\left(\frac{-1}{A_2}\right)B_2)+(-\mu\alpha_2)(\left(\frac{-1}{A_2}\right)A_2C_1)=0$
\item $\left(\frac{-1}{A_2}\right)A_2C_1\equiv -1$ (mod $4$)
\item $-\mu\left(\frac{-1}{A_2}\right)A_1C_2\equiv -1$ (mod $4$)
\item $(\alpha_1,B_1,\left(\frac{-1}{A_2}\right)A_2C_1)=1$
\item $(-\mu\alpha_2,-\mu\left(\frac{-1}{A_2}\right)B_2,-\mu\left(\frac{-1}{A_2}\right)A_1C_2)=1$
\end{enumerate}
\end{comment}

\begin{tabular}{ll}
(1) $\alpha_1(-\mu\left(\frac{-1}{A_2}\right)A_1C_2)+4(B_1)(-\mu\left(\frac{-1}{A_2}\right)$&\hspace{-.3cm}$B_2)+(-\mu\alpha_2)(\left(\frac{-1}{A_2}\right)A_2C_1)=0$,\\
(2) $\left(\frac{-1}{A_2}\right)A_2C_1\equiv -1$ (mod $4$),&
(3) $-\mu\left(\frac{-1}{A_2}\right)A_1C_2\equiv -1$ (mod $4$),\\
(4) $(\alpha_1,B_1,\left(\frac{-1}{A_2}\right)A_2C_1)=1$,&
(5) $(-\mu\alpha_2,-\mu\left(\frac{-1}{A_2}\right)B_2,-\mu\left(\frac{-1}{A_2}\right)A_1C_2)=1.$
\end{tabular}

We will omit these details. Putting everything together we find that the pair of 6-tuples from line (\ref{6tuples}) is an element of $\Gamma_\infty\backslash \mathbb{S}(A_1,\mu A_2) \times\Gamma_\infty\backslash \mathbb{S}(\alpha_1,-\mu\alpha_2).$  Now we can apply the canonical map from $\Gamma_{\infty}\backslash\mathbb{S}(\cdot,\cdot)$ to $\Gamma_{\infty}\backslash\mathbb{S}(\cdot,\cdot)/\Gamma_{\infty}$. Let
\[\phi:\Gamma_\infty\backslash \mathbb{S}(A_1\alpha_1,A_2\alpha_2)/\Gamma_\infty \rightarrow \Gamma_\infty\backslash \mathbb{S}(A_1,\mu A_2)/\Gamma_\infty \times\Gamma_\infty\backslash \mathbb{S}(\alpha_1,-\mu\alpha_2)/\Gamma_\infty\]
be the map constructed above.

To establish the bijection we construct the inverse map $\psi$.  Consider
 \[((4A_1,4b_1,c_1,\mu 4A_2,4b_2,c_2),(4\alpha_1,4\beta_1,\gamma_1,-\mu4\alpha_2,4\beta_2,\gamma_2))\in S(A_1,\mu A_2)\times S(\alpha_1,-\mu\alpha_2).\]
Using the bijection described in Proposition \ref{DoubleCosetProp}, we will identify this element of $S(A_1,\mu A_2)\times S(\alpha_1,-\mu\alpha_2)$ with its corresponding element in $\Gamma_{\infty}\backslash\mathbb{S}(A_1,\mu A_2)/\Gamma_{\infty}\times \Gamma_{\infty}\backslash \mathbb{S}(\alpha_1,-\mu\alpha_2)/\Gamma_{\infty}$. 
Let 
\begin{multline}\psi(((4A_1,4b_1,c_1,\mu 4A_2,4b_2,c_2),(4\alpha_1,4\beta_1,\gamma_1,-\mu4\alpha_2,4\beta_2,\gamma_2))\\=(4A_1\alpha_1,4B_1,C_1,4A_2\alpha_2,4B_2,C_2).\label{psidef}\end{multline}
We define $B_i,C_i$ presently. Define $B_1$ to be the integer such that $B_1\equiv b_1$ (mod $A_1$), $B_1\equiv \beta_1$ (mod $\alpha_1$), and such that $B_1/A_1\alpha_1\in[0,1).$ Define $B_2$ to be the integer such that $B_2\equiv b_2$ (mod $A_2$), $B_2\equiv -\mu\left(\frac{-1}{A_2}\right)\beta_2$ (mod $\alpha_2$), and such that $B_2/A_2\alpha_2\in[0,1).$  

The definition of $C_2$ will be a bit more involved. We begin by rewriting the congruence conditions of the previous paragraph as equalities; specifically, $b_1=B_1-xA_1,$ $b_2=B_2+y\mu A_2,$ $\beta_1=B_1-x^\prime\alpha_1,$ and $\beta_2=-\mu\left(\frac{-1}{A_2}\right)B_2+y^\prime(-\mu)\alpha_2$, where $x,y,x^\prime,$ and $y^\prime\in\mathbb{Z}$. Let $C_2$ be defined to be the unique integer satisfying $C_2\equiv \gamma^{-1}(c_2-4B_2x)$ (mod $4A_2$), $C_2\equiv -\mu\left(\frac{-1}{A_2}\right)A_1^{-1}(\gamma_2-4(-\mu\left(\frac{-1}{A_2}\right))B_2x^\prime)$ (mod $4\alpha_2$), and $C_2/4A_2\alpha_2\in[0,1)$. For later we record two more equalities; let $z,z^\prime\in\mathbb{Z}$ such that $c_2=\gamma C_2+4B_2x-4\mu A_2z$ and $\gamma_2 = -\mu\left(\frac{-1}{A_2}\right)A_1C_2-4\mu\left(\frac{-1}{A_2}\right)B_2x^\prime+4\mu\alpha_2z^\prime$. 

Finally define $C_1=\frac{-A_1\alpha_1C_2-4B_1B_2}{A_2\alpha_2}$. 

To show that $(4A_1\alpha_1,4B_1,C_1,4A_2\alpha_2,4B_2,C_2)\in S(A_1\alpha_1,A_2\alpha_2)$, we must verify the following six conditions:
\begin{enumerate}
\item $A_1\alpha_1C_2+4B_1B_2+A_2\alpha_2C_1=0$
\item $C_2\equiv -1$ (mod $4$)
\item $(A_2\alpha_2,B_2,C_2)=1$
\item $C_1\in\mathbb{Z}$
\item $C_1\equiv -1$ (mod $4$)
\item $(A_1\alpha_1,B_1,C_1) = 1$

\end{enumerate}

For this list we will provide a proof of item (6) and make a comment about items (4) and (5).

Item (4) can be proved by analyzing the numerator of $C_1=\frac{-A_1\alpha_1C_2-4B_1B_2}{A_2\alpha_2}$ modulo $A_2$ and $\alpha_2$ separately; the (perhaps unintuitive) congruence conditions defining $C_2$ play an important role in this calculation. 

To continue we require two additional identities involving $C_1$.  The first is 
\begin{equation}
C_1 = \frac{\gamma_1+4B_1y^\prime+4\alpha_1(-x^\prime y^\prime-z^\prime)}{\left(\frac{-1}{A_2}\right)A_2}.\label{C1first}
\end{equation}  The second is 
\begin{equation}
C_1=\frac{c_1+B_1y-A_1(xy+z)-\mu\ell A_1C_2}{\mu\alpha_2}.\label{C1second}
\end{equation}
We postpone the proof of these identities for the moment. 

Item (5) can be proved using equation (\ref{C1first}).

To prove item (6), we begin by observing that  $(A_1\alpha_1,B_1,C_1)=(A_1,B_1,C_1)(\alpha_1,B_1,C_1)$. To see that $(\alpha_1,B_1,C_1) = 1$ recall that $B_1=\beta_1+x^\prime\alpha_1$, $C_1= \frac{\gamma_1+4B_1y^\prime+4\alpha_1(-x^\prime y^\prime-z^\prime)}{\left(\frac{-1}{A_2}\right)A_2}$, and $(A_2,\alpha_1)=1$. These results imply that 
\begin{equation*}(\alpha_1,B_1,C_1)
=(\alpha_1,\beta_1+x^\prime \alpha_1,\frac{\gamma_1+4B_1y^\prime+4\alpha_1(-x^\prime y^\prime-z^\prime)}{\left(\frac{-1}{A_2}\right)A_2})
=(\alpha_1,\beta_1,\gamma_1)
=1.
\end{equation*}

It remains to show that $(A_1,B_1,C_1) =1.$ This is similar to the previous calculation, but we must use the second identity for $C_1$.

%%%NOTE: The extra kA_1C_2 appearing in the second formula for C_1 has to do with the gamme that is needed in the congruence condition used to define C_2.
%%%

Next we will prove the two equalities involving $C_1$. First we address equation (\ref{C1first}).
 
 Since $\beta_1=B_1-x^\prime\alpha_1$,
\begin{align}
A_2\alpha_2\text{RHS}((\ref{C1first}))\label{spltempRHS}=&\left(\frac{-1}{A_2}\right)(\gamma_1\alpha_2+4\beta_1\alpha_2y^\prime-4\alpha_1\alpha_2z^\prime).
\end{align}

Before addressing the left hand side of the equality, recall that $-\mu\left(\frac{-1}{A_2}\right)A_1C_2=\gamma_2+4\mu\left(\frac{-1}{A_2}\right)B_2x^\prime-4\mu\alpha_2z^\prime$. Thus,
\begin{align}
A_2\alpha_2\text{LHS}((\ref{C1first}))
\label{spltemp4}=&\alpha_1\mu\left(\frac{-1}{A_2}\right)\gamma_2-\alpha_1\left(\frac{-1}{A_2}\right)4\alpha_2z^\prime-4B_2\beta_1\\
\label{spltempLHS}=&\mu\left(\frac{-1}{A_2}\right)(\alpha_1\gamma_2+4\beta_1\beta_2)-4\left(\frac{-1}{A_2}\right)\alpha_1\alpha_2z^\prime+4\left(\frac{-1}{A_2}\right)y^\prime\alpha_2\beta_1.
\end{align}

Next we apply the identity $\alpha_1\gamma_2+4\beta_1\beta_2-\mu\alpha_2\gamma_1=0$ to see that (\ref{spltempRHS}) is equal to (\ref{spltempLHS}). So we have established equation (\ref{C1first}).

The proof of equation (\ref{C1second}) is similar and will be omitted.

Finally we must show that the maps $\phi$ and $\psi$ are inverses. We begin by computing $\psi\circ\phi$. Let $(4A_1\alpha_1,4B_1,C_1,4A_2\alpha_2,4B_2,C_2)\in S(A_1\alpha_1,A_2\alpha_2)$.  The image of this point under the map $\phi$ is  
\begin{multline}((4A_1,4B_1,\mathcal{C}_1,\mu 4A_2,4B_2,\gamma C_2),\\(4\alpha_1,4B_1,\left(\frac{-1}{A_2}\right) A_2C_1,-\mu4\alpha_2,-\left(\frac{-1}{A_2}\right)\mu 4B_2,-\mu \left(\frac{-1}{A_2}\right)A_1C_2)),\end{multline}
where $\mathcal{C}_1=\frac{-A_1\gamma C_2-4B_1B_2}{\mu A_2}$.
Now let $x,x^\prime,y,y^\prime,z,z^\prime\in\mathbb{Z}$ such that 
\begin{multline*}((4A_1,4B_1-4A_1x,\mathcal{C}_1-4B_1y+4A_1(xy-z),\mu 4A_2,4B_2+\mu 4A_2y,\gamma C_2+4B_2x+\mu 4A_2z),\\(4\alpha_1,4B_1-4\alpha_1x^\prime,\left(\frac{-1}{A_2}\right) A_2C_1-4B_1y^\prime+4\alpha_1(x^\prime y^\prime-z^\prime),\\-\mu4\alpha_2,-\left(\frac{-1}{A_2}\right)\mu 4B_2-\mu 4\alpha_2y^\prime,-\mu \left(\frac{-1}{A_2}\right)A_1C_2-\left(\frac{-1}{A_2}\right)\mu 4B_2x^\prime-\mu4\alpha_2z^\prime))\end{multline*}
is an element of $S(A_1,\mu A_2)\times S(\alpha_1,-\mu\alpha_2)$. Now the image of this element under the map $\psi$ is 
\begin{equation*}
(4A_1\alpha_1,4B_1^*,C_1^*,4A_2\alpha_2,4B_2^*,C_2^*),
\end{equation*}
where the definitions of $B_1^*$ ,$C_1^*$, $B_2^*$, and $C_2^*$ can be found in the two paragraphs following line (\ref{psidef}). We claim that $(4A_1\alpha_1,4B_1^*,C_1^*,4A_2\alpha_2,4B_2^*,C_2^*)=(4A_1\alpha_1,4B_1,C_1,4A_2\alpha_2,4B_2,C_2)$. This equality is a consequence of the Chinese remainder theorem. Thus $\psi\circ\phi = \operatorname{id}$.

The last thing that we must check is that $\phi\circ\psi=\operatorname{id}$. Consider
 \[((4A_1,4b_1,c_1,\mu 4A_2,4b_2,c_2),(4\alpha_1,4\beta_1,\gamma_1,-\mu4\alpha_2,4\beta_2,\gamma_2))\in S(A_1,\mu A_2)\times S(\alpha_1,-\mu\alpha_2).\]
Let 
\begin{multline}\psi(((4A_1,4b_1,c_1,\mu 4A_2,4b_2,c_2),(4\alpha_1,4\beta_1,\gamma_1,-\mu4\alpha_2,4\beta_2,\gamma_2))\\=(4A_1\alpha_1,4B_1,C_1,4A_2\alpha_2,4B_2,C_2),\end{multline}
where $B_1,C_1,B_2,C_2$ are defined as in the two paragraphs following line (\ref{psidef}). Let the image of this element under the map $\phi$ be
\begin{multline*}((4A_1,4B_1-4A_1x,\mathcal{C}_1-4B_1y+4A_1(xy-z),\mu 4A_2,4B_2+\mu 4A_2y,\gamma C_2+4B_2x+\mu 4A_2z),\\(4\alpha_1,4B_1-4\alpha_1x^\prime,\left(\frac{-1}{A_2}\right) A_2C_1-4B_1y^\prime+4\alpha_1(x^\prime y^\prime-z^\prime),\\-\mu4\alpha_2,-\left(\frac{-1}{A_2}\right)\mu 4B_2+\mu \alpha_2y^\prime,-\mu \left(\frac{-1}{A_2}\right)A_1C_2-\left(\frac{-1}{A_2}\right)\mu 4B_2x^\prime-\mu4\alpha_2z^\prime)),\end{multline*}
where $x,x^\prime,y,y^\prime,z,z^\prime\in\mathbb{Z}$ such that this pair of 6-tuples is an element of $S(A_1,\mu A_2)\times S(\alpha_1,-\mu\alpha_2)$. Again we use the Chinese remainder theorem to see that $\phi\circ\psi=\operatorname{id}$.\EndProof

\subsection{The Splitting in Terms of Block Parameters} \label{sec:BlockSplit}

This brief section recalls Miller's formula \cite{M06} for the splitting in terms of block parameters. Suppose that $\gamma\in \Gamma$ and $\gamma = n
\left( \begin{smallmatrix}
a_1 & b_1 & 0\\
c_1 & d_1 & 0\\
0 & 0 & 1 \end{smallmatrix} \right)
\left( \begin{smallmatrix}
a_2 & 0 & b_2\\
0 & 1 & 0\\
c_2 & 0 & d_2 \end{smallmatrix} \right)\
\left( \begin{smallmatrix}
1 & 0 & 0\\
0 & a_3 & b_3\\
0 & c_3 & d_3 \end{smallmatrix} \right)$, where each 
$\left(\begin{smallmatrix}
a_i & b_i\\
c_i & d_i\end{smallmatrix}\right)\in \Gamma_1(4)\subseteq \text{SL}(2,\mathbb{Z}).$
Then 
\begin{equation} \label{blocksplitting} s(\gamma) = \left(\frac{c_1}{d_1}\right)\left(\frac{c_2}{d_2}\right)\left(\frac{c_3}{d_3}\right)s_{NA}(\gamma),\end{equation}
where 
\begin{equation} \label{nonarithsplitting} s_{NA}(\gamma)=\begin{cases}
(c_1,d_1) &,\hspace{.5cm} c_1,c_2\neq 0,\, c_3=0 \\
(c_1c_2(c_1c_3-d_1c_2a_3),c_1a_3)(a_3,c_2c_3) &,\hspace{.5cm} c_1,c_2,c_3, c_1c_3-d_1c_2a_3\neq0\\
(c_2a_3,c_1a_3)(a_3,c_2c_3) &,\hspace{.5cm}c_1,c_2,c_3\neq0,\, c_1c_3-d_1c_2a_3=0\\
(a_3,c_2c_3) &,\hspace{.5cm}c_1=0,\,c_2,c_3\neq0\\
1 &,\hspace{.5cm}\text{otherwise}.\end{cases}\end{equation}

\noindent First, note that this formula shows that $s$ is left $\Gamma_{\infty}$-invariant. Second, note that $s$ is not a group homomorphism, but rather that $s(\gamma_1\gamma_2)=s(\gamma_1)s(\gamma_2)\sigma(\gamma_1,\gamma_2)$.

%%%%%%%%%%%%%%%%%%%%%%%%%%%%%%%%%%%%
%%%%%%%Symmetries of the Splitting
%%%%%%%%%%%%%%%%%%%%%%%%%%%%%%%%%%%%

\subsection{Symmetries of the Splitting} \label{sec:SymSpl}

This section describes how the map $s$ is affected by some of the symmetries described in section \ref{sec:PluckerSym}. 

\subsection{Conjugation} \label{ssec:Conj}

The simplest automorphisms that preserve $\Gamma_{\infty}$ are given by conjugation by elements of $\Gamma_{\infty}$. This action will leave the splitting unchanged.

\begin{proposition}\label{doublecosetsym} Let $\gamma\in \Gamma_{1}(4)$ and $n\in \Gamma_{\infty}$. Then $s(n\gamma n^{-1}) = s(\gamma)$. In fact, $s(n\gamma)=s(\gamma n)=s(\gamma)$.
\end{proposition}

\noindent\textbf{Proof:} We will show that $s(\gamma n)=s(\gamma)$. The case of $s(n\gamma) = s(\gamma)$ is similar. Now by the definition of the splitting $s(\gamma n) = s(\gamma)s(n)\sigma(\gamma,n)$.  By results from  \cite{M06} and  \cite{BLS99}, $s(n) = 1$ and $\sigma(\gamma,n) = 1$.  Therefore $s(\gamma n) = s(\gamma)$. \EndProof

\noindent Proposition \ref{doublecosetsym} shows that $s$ is well defined on the double coset space $\Gamma_{\infty}\backslash\Gamma_{1}(4)/\Gamma_{\infty}$. 

Next consider conjugation by the elements $S_2=\text{t}(1,-1,1)$ and $S_3 = \text{t}(1,1,-1)$. The following proposition is not exhaustive; it only includes symmetries that will be used in the sequel.

\noindent\begin{proposition} \label{SignSym} Let $\gamma\in\Gamma$. Then
\[s(S_3\gamma S_3) = \begin{cases}
(c_1a_3,-1)s(\gamma) &,\, c_1,c_2,c_3\neq0,A_2=0\\
s(\gamma) &,\text{ otherwise}
\end{cases}\]
If $A_1,A_2\neq 0.$ Then $s(S_2\gamma S_2)=-\text{sign}(A_1A_2)s(\gamma)$.
\end{proposition}

The proof of this result is similar to the proof in the next subsection and will be omitted. 

%%%%%%%%%%%%%%%%%
%%%%%%%%%%%%%%%%%
%%%%%%%%%%%%%%%%%

\subsection{Cartan Involution Composed with the Long Element} \label{ssec:CartLong}

\begin{proposition} \label{CartSym} Let $\gamma\in \Gamma_1(4)$ with Pl{\"u}cker coordinates $(4A_1,4B_1,C_1,4A_2,4B_2,C_2).$  Consider the involution $\phi:\gamma\mapsto w_{\ell} \gamma^{-t}w_{\ell}^{-1}$. When $A_1$ and $A_2$ are not equal to $0$ 
\[s(\phi(\gamma)) = (-A_1,-A_2)s(\gamma).\]
When $A_1,B_2\neq0$ and $A_2=0$,
\[s(\phi(\gamma))=(-A_1,B_2)s(\gamma).\]
\end{proposition}
\noindent A similar result holds on the other cells, but these identities will not be needed.

\noindent\textbf{Proof:} By equation (\ref{blocksplitting}), if $\gamma = n\gamma_1\gamma_2\gamma_3$, then $s(\gamma)=s(\gamma_1)s(\gamma_2)s(\gamma_3)\sigma(\gamma_1,\gamma_2\gamma_3)\sigma(\gamma_2,\gamma_3)$, and $s(\phi(\gamma)) = s(\phi(\gamma_1))s(\phi(\gamma_2))s(\phi(\gamma_3))\sigma(\phi(\gamma_1),\phi(\gamma_2)\phi(\gamma_3))\sigma(\phi(\gamma_2),\phi(\gamma_3))$. Direct calculation shows that  
\[\begin{array}{ccccc}
\phi(\gamma_1)&=&\phi\left(\left( \begin{smallmatrix}
a_1 & b_1 & 0\\
c_1 & d_1 & 0\\
0 & 0 & 1 \end{smallmatrix} \right)\right)&=&\left( \begin{smallmatrix}
1 & 0 & 0\\
0 & a_1 & b_1\\
0 & c_1 & d_1 \end{smallmatrix} \right),\\
\phi(\gamma_2)&=&\phi\left(\left( \begin{smallmatrix}
a_2 & 0 & b_2\\
0 & 1 & 0\\
c_2 & 0 & d_2 \end{smallmatrix} \right)\right)&=&
\left( \begin{smallmatrix}
a_2 & 0 & -b_2\\
0 & 1 & 0\\
-c_2 & 0 & d_2 \end{smallmatrix} \right),\\
\phi(\gamma_3)&=&\phi\left(\left(\begin{smallmatrix}
1 & 0 & 0\\
0 & a_3 & b_3\\
0 & c_3 & d_3 \end{smallmatrix}\right)\right)&=&\left( \begin{smallmatrix}
a_3 & b_3 & 0\\
c_3 & d_3 & 0\\
0 & 0 & 1 \end{smallmatrix} \right).\end{array}\]
Thus $s(\phi(\gamma_1)) = s(\gamma_1)$, $s(\phi(\gamma_3))=s(\gamma_3)$, and $s(\phi(\gamma_2)) = s(\gamma_2)$. 

It remains to compute $\sigma(\gamma_1,\gamma_2\gamma_3)\sigma(\gamma_2,\gamma_3)\sigma(\phi(\gamma_1),\phi(\gamma_2)\phi(\gamma_3))\sigma(\phi(\gamma_2),\phi(\gamma_3))$.  First note that Proposition 4.2 in \cite{M06} proves that  

\begin{equation}\label{siggam2gam3}
\sigma(\gamma_2,\gamma_3)=
\begin{cases}
(a_3,c_2c_3), &\,c_2,c_3\neq0\\
1, &\,\text{otherwise.}
\end{cases}
\end{equation}
The details of this computation will be omitted as they resemble what is about to follow.

Next we compute $\sigma(\phi(\gamma_2),\phi(\gamma_3))$. If $c_3=0$, then $\gamma_3\in N$ and $\sigma(\phi(\gamma_2),\phi(\gamma_3))=1$. Assume that $c_2\neq0$ and $c_3\neq 0$. Using the identity $\sigma(n_1g_1n_2,g_2n_3)=\sigma(g_1,n_2g_2)$ we can show that
\begin{align}
\nonumber\sigma(\phi(\gamma_2),\phi(\gamma_3))=&\sigma\left 
(\left( \begin{smallmatrix}
 &  & c_2^{-1}\\
 & 1 & \\
-c_2 &  & d_2 \end{smallmatrix} \right),
\left( \begin{smallmatrix}
a_3 &  & \\
c_3 & a_3^{-1} & \\
 &  & 1 \end{smallmatrix} \right)\right)\\
\nonumber=&\sigma\left(
\left( \begin{smallmatrix}
c_2^{-1} &  & \\
 & -1 & \\
 &  & -c_2 \end{smallmatrix} \right)w_{\alpha_1}w_{\alpha_2}w_{\alpha_1},
\left( \begin{smallmatrix}
1 &  & \frac{-d_2}{c_2}\\
 & 1 & \\
 &  & 1 \end{smallmatrix} \right)
\left( \begin{smallmatrix}
a_3 &  & \\
c_3 & a_3^{-1}& \\
 &  & 1 \end{smallmatrix} \right)\right)\\
=&\sigma\left(
\left(\begin{smallmatrix}
c_2^{-1} &  & \\
 & -1 & \\
 &  & -c_2 \end{smallmatrix} \right)w_{\alpha_1}w_{\alpha_2}w_{\alpha_1},
\left( \begin{smallmatrix}
a_3 &  & \\
c_3 & a_3^{-1}& \label{eqn1} \\
 &  & 1 \end{smallmatrix} \right)\right)\end{align}
Now using equation (\ref{BLSdef}) we see that 
 \begin{align}
\nonumber(\ref{eqn1})
\label{eqn2}=&(a_3,-1)(a_3c_3,-c_3)= (a_3,-1)(a_3,-c_3)(c_3,-c_3) = (a_3,c_3)
\end{align}

If either $c_2$ or $c_3$ is equal to $0$, then $\sigma(\phi(\gamma_2),\phi(\gamma_3)) = 1.$  Thus
\begin{equation*}\sigma(\phi(\gamma_2),\phi(\gamma_3))\sigma(\gamma_2,\gamma_3)=\bigg\{\begin{smallmatrix} 
(a_3,c_2), & c_2,c_3\neq 0\\
1, & \text{otherwise}.\end{smallmatrix}
\end{equation*}
However, if $c_3=0$ it follows that $a_3=1$ so in fact we have 
\begin{equation}
\sigma(\phi(\gamma_2),\phi(\gamma_3))\sigma(\gamma_2,\gamma_3)=(a_3,c_2)\label{tempsplitcocycle1}
\end{equation}
%At this point note that when $c_1=0$, then $a_3 = -A_2.$  Thus this computation agrees with the statement of the proposition.  
It remains to consider $\sigma(\gamma_1,\gamma_2\gamma_3)$ and $\sigma(\phi(\gamma_1),\phi(\gamma_2)\phi(\gamma_3))$. In Proposition 4.2 in \cite{M06}, Miller proves that
\begin{equation*}
\sigma(\gamma_1,\gamma_2\gamma_3)=
\begin{cases}
(c_1c_2(-A_2),c_1a_3),&\,c_1,c_2,A_2\neq 0\\
(c_2a_3,c_1a_3),&\,c_1,c_2\neq 0, A_2=0\\
1,&\,\text{otherwise.}
\end{cases}
\end{equation*}

Now we will compute $\sigma(\phi(\gamma_1),\phi(\gamma_2)\phi(\gamma_3))$. Note that if $c_1=0$, then this 2-cocycle is equal to 1. Thus assume that $c_1\neq0$. As $a_j\equiv 1 \text{ } (4)$, thus $a_j\neq 0$.
In the following computation let 
\[\begin{array}{lll}
\alpha&=&\text{the first nonzero quantity among } c_2, \frac{-A_2}{c_1a_2}, \frac{1}{a_2a_3},\\
\beta&=&\text{the first nonzero quantity among} \frac{-A_2}{c_1}, \frac{1}{a_3},\\
\delta&=&\text{the first nonzero quantity among} -c_2a_3, \frac{1}{a_2},\\
h&=&\left( \begin{smallmatrix}
a_2a_3& & \\
\frac{-A_2}{c_1}& a_3^{-1} &\\
-c_2a_3& &a_2^{-1} \end{smallmatrix}\right).\end{array}\]
Then using the identity $\sigma(n_1g_1n_2,g_2n_3)=\sigma(g_1,n_2g_2)$ it follows that
\begin{align}
\nonumber\sigma(\phi(\gamma_1),\phi(\gamma_2)\phi(\gamma_3))
=&\sigma\left(\left(\begin{smallmatrix}
1& &\\
& &-c_1^{-1}\\
& c_1&d_1\end{smallmatrix}\right),
\left( \begin{smallmatrix}
a_2a_3& a_2b_3& -b_2\\
c_3& d_3& \\
-c_2a_3& -c_2b_3& d_2\end{smallmatrix}\right)\right)\\
=&\sigma\left(\left(\begin{smallmatrix}
1& &\\
& &-c_1^{-1}\\
& c_1&d_1\end{smallmatrix}\right),
\left( \begin{smallmatrix}
a_2a_3& & \\
c_3& a_3^{-1} &\\
-c_2a_3& &a_2^{-1} \end{smallmatrix}\right)
\right).\label{place1}
\end{align}
Now factor the matrix in the first entry using the SL$(2,\mathbb{R})$ Bruhat decomposition on the lower right $2\times2$ block and apply the identity $\sigma(n_1g_1n_2,g_2n_3)=\sigma(g_1,n_2g_2)$ to see that
\begin{align}
\nonumber(\ref{place1})=&\sigma\left(
\left(\begin{smallmatrix}
1& &\\
& c_1^{-1}& \\
& & c_1\end{smallmatrix}\right)
\left(\begin{smallmatrix}
1& & \\
& &-1 \\
& 1& \end{smallmatrix}\right)
\left(\begin{smallmatrix}
1& & \\
& 1&\frac{d_1}{c_1} \\
& & 1\end{smallmatrix}\right)
,\left( \begin{smallmatrix}
a_2a_3& & \\
c_3& a_3^{-1} &\\
-c_2a_3& &a_2^{-1} \end{smallmatrix}\right)\right)\\
=&\sigma\left(
\left(\begin{smallmatrix}
1& &\\
& c_1^{-1}& \\
& & c_1\end{smallmatrix}\right)
\left(\begin{smallmatrix}
1& & \\
& &-1 \\
& 1& \end{smallmatrix}\right),
\left( \begin{smallmatrix}
a_2a_3& & \\
\frac{-A_2}{c_1}& a_3^{-1} &\\
-c_2a_3& &a_2^{-1} \end{smallmatrix}\right)
\right).\label{eqn3}\end{align}
Now equation (\ref{BLSdef}) shows that 
\begin{align}
\nonumber(\ref{eqn3})=&(c_1,\beta)(\beta\delta,-\delta)=(-c_1\delta,\beta).\label{eqn4}
\end{align}
Thus, 
\begin{equation}\sigma(\phi(\gamma_1),\phi(\gamma_2)\phi(\gamma_3))\sigma(\gamma_1,\gamma_2\gamma_3)= (\alpha,\beta c_1\delta).\label{tempsplitcocycle4}
\end{equation}
If we combine lines (\ref{tempsplitcocycle1}) and (\ref{tempsplitcocycle4}) and use the fact that $A_1\neq 0$ we get
\begin{multline}
\sigma(\phi(\gamma_2),\phi(\gamma_3))\sigma(\gamma_2,\gamma_3)\sigma(\phi(\gamma_1),\phi(\gamma_2)\phi(\gamma_3))\sigma(\gamma_1,\gamma_2\gamma_3)\\=(c_2,a_3)(c_2,\beta c_1(-c_2a_3)) = (c_2,\beta c_1).\label{tempcocycle5}
\end{multline}
Remember that this equation is only valid if $c_1\neq0$. If $c_1=0$, then as mentioned above $\sigma(\phi(\gamma_1),\phi(\gamma_2)\phi(\gamma_3))\sigma(\gamma_1,\gamma_2\gamma_3)=1$ and so we get 
\begin{equation}
\sigma(\phi(\gamma_2),\phi(\gamma_3))\sigma(\gamma_2,\gamma_3)\sigma(\phi(\gamma_1),\phi(\gamma_2)\phi(\gamma_3))\sigma(\gamma_1,\gamma_2\gamma_3)=(c_2,a_3).\label{tempcocycle6}
\end{equation}

If $A_1,A_2,c_1\neq 0$, then we have
 \[(\ref{tempcocycle5})= (c_2,\beta c_1) = (-A_1,-A_2).\]
 
If $A_1,A_2\neq0, c_1=0$, then $d_1=1$ and $A_2=-A_1a_3$. In this case we have
\[(\ref{tempcocycle6})=(c_2,a_3)=(-A_1,-A_1A_2) = (-A_1,-A_2).\]
 
 If $A_2=0$, Then $c_1,c_3\neq0$ and $c_1c_3=d_1c_2a_3$. If additionally, $A_1,B_2\neq0$, then we have
 \begin{multline*}(\ref{tempcocycle5}) = (c_2,a_3 c_1)=(c_2,c_2c_3d_1)=(c_2,-c_3d_1)\\
= (-A_1,-C_2B_1)=(-A_1,-C_2(4B_1B_2)B_2)=(-A_1,C_2^2A_1B_2)=(-A_1,B_2).
 \end{multline*}\EndProof
 
%%%%%%%%%%%%%%%%%%%%%%%%
%%%%%%%%%Splitting
%%%%%%%%%%%%%%%%%%%%%%%%

\subsection{The Splitting in Terms of Pl\"{u}cker Coordinates} \label{sec:PluckSplit}

\begin{theorem} \label{theorem:PluckSplit} Let $\gamma\in\Gamma_{1}(4)$ with Pl\"{u}cker coordinates $(4A_1,4B_1,C_1,4A_2,4B_2,C_2)$ such that $A_1>0$, and $A_2/(A_1,A_2)\equiv 1\,(\text{mod } 2)$. Let $D=(A_1,A_2)$, $D_1 = (D,B_1)$, $D_2=D/D_1$, and let $\epsilon = \left(\frac{-1}{-B_1/D_1}\right)$. Then
\begin{equation} \label{plucksplitting}
s(\gamma)=\left(\frac{-\epsilon}{-A_1A_2}\right)\left(\frac{A_1/D}{A_2/D}\right)
\left(\frac{B_1/D_1}{A_1/D}\right)\left(\frac{4B_2/D_2}{\text{sign}(A_2)A_2/D}\right)\left(\frac{D_1}{C_1}\right)\left(\frac{D_2}{C_2}\right).\end{equation}
\end{theorem}
The proof of this theorem occupies the next two subsections. Subsection \ref{sec:SplitPlucker2} contains the reduction step and Subsection \ref{sec:SplitPlucker3} contains the computation on reduced elements. 

A few remarks are in order regarding the formula for the splitting. First, this formula can be used to derive the value of the splitting on any input using the symmetries of Subsection \ref{sec:SymSpl}. This point will be discussed more throughly in Subsection \ref{sec:SplitPlucker3}. Second, if $-1$ is assumed to be a square then this formula essentially reduces to the formula of Brubaker-Bump-Friedberg-Hoffstein \cite{BBFH07}. Third, $\left(\frac{-\epsilon}{-A_1A_2}\right)=1$ when $A_1$ and $A_2$ are odd.

%%%%%%%%%%%%%%%%%%%%%%%%
%%%%%%%%%Reduction
%%%%%%%%%%%%%%%%%%%%%%%%

\subsection{The Splitting: The Reduction} \label{sec:SplitPlucker2}

\begin{proposition} \label{splitreduction} Let $\gamma\in \Gamma_1(4)$ with Pl\"{u}cker Coordinates $(4A_1,4B_1,C_1,4A_2,4B_2,C_2)$, such that $A_1,A_2\neq0$.  Suppose that $D$ divides $(A_1,A_2)$.  Let $D_1 = (D,B_1)$ and let $D_2=D/D_1.$  Suppose that $D_2$ divides $B_2$. Let $ S=\left( \begin{smallmatrix}
1 & 0 & 0 \\
0 & D_2^{-1} & 0 \\
0 & 0 & D^{-1} \end{smallmatrix} \right).$
Then $S\gamma S^{-1}\in \Gamma_{1}(4)$ with Pl\"{u}cker Coordinates $(4A_1/D,4B_1/D_1,C_1,4A_2/D,4B_2/D_2,C_2)$ and
\begin{equation}
s(\gamma) = s(S\gamma S^{-1})\left(\frac{D_1}{C_1}\right)\left(\frac{D_2}{C_2}\right).
\end{equation}
\end{proposition}

\noindent\textbf{Proof:} As the Pl\"{u}cker coordinates satisfy the previously mentioned divisibility conditions, $\gamma$ is of the form $\left( \begin{smallmatrix}
a_{11} & a_{12} & a_{13} \\
D_2a_{21} & a_{22} & a_{23} \\
Da_{31} & D_1a_{32} & a_{33} \end{smallmatrix} \right)$. If $\gamma$ is factored into blocks, then \newline$\gamma = n
\left( \begin{smallmatrix}
a_1 & b_1 & 0\\
D_2c_1 & d_1 & 0\\
0 & 0 & 1 \end{smallmatrix} \right)
\left( \begin{smallmatrix}
a_2 & 0 & b_2\\
0 & 1 & 0\\
Dc_2 & 0 & d_2 \end{smallmatrix} \right)\
\left( \begin{smallmatrix}
1 & 0 & 0\\
0 & a_3 & b_3\\
0 & D_1c_3 & d_3 \end{smallmatrix} \right)$.  
%%%%%%%%%%%%%%%%%%%%%%%%%%%%%%%%%%%%%%
%%%%%%%%%That the factorization works as described should be proved.
%%%%%%%%%%%%%%%%%%%%%%%%%%%%%%%%%%%%%%
Using Proposition \ref{prop:PluckerSym} we see that $\gamma\in \mathbb{S}(A_1D,A_2D)$ is mapped into $\mathbb{S}(A_1,A_2)$ via conjugation by the matrix $S=\left( \begin{smallmatrix}
1 & 0 & 0 \\
0 & D_2^{-1} & 0 \\
0 & 0 & D^{-1} \end{smallmatrix} \right)$. Explicitly, $S\gamma_1S^{-1}=\left( \begin{smallmatrix}
a_1 & D_2b_1 & 0 \\
c_1 & d_1 & 0 \\
0 & 0 & 1 \end{smallmatrix} \right)$, 
$S\gamma_2S^{-1}=\left( \begin{smallmatrix}
a_2 & 0 & Db_2 \\
0 & 1 & 0 \\
c_2 & 0 & d_2 \end{smallmatrix} \right)$, and
$S\gamma_3S^{-1}=\left( \begin{smallmatrix}
1 & 0 & 0 \\
0 & a_3 & D_1b_3 \\
0 & c_3 & d_3 \end{smallmatrix} \right)$.
 
Using equation (\ref{blocksplitting}),
\begin{align*}
s(S\gamma S^{-1})=&\left(\frac{c_1D_2}{d_1}\right)\left(\frac{c_2D}{d_2}\right)\left(\frac{c_3D_1}{d_3}\right)\left(\frac{D_2}{d_1}\right)\left(\frac{D}{d_2}\right)\left(\frac{D_1}{d_3}\right)\\
&\hspace{2cm}\times\sigma(S\gamma_1S^{-1},S\gamma_2\gamma_3S^{-1})\sigma(S\gamma_2S^{-1},S\gamma_3S^{-1})\\
=&s(\gamma)\left(\frac{D_2}{-C_2/d_2}\right)\left(\frac{D}{d_2}\right)\left(\frac{D_1}{-C_1/d_2}\right)=s(\gamma)\left(\frac{D_1}{-C_1}\right)\left(\frac{D_2}{-C_2}\right).
\end{align*}
The third equality follows from formula (\ref{nonarithsplitting}) paired with the fact that conjugation by S only scales the block parameters by positive numbers leaving the non-arithmetic factor unchanged. It follows that $s(\gamma) = s(S\gamma S^{-1})\left(\frac{D_1}{C_1}\right)\left(\frac{D_2}{C_2}\right)$.\EndProof
\begin{comment}
Using the computation from section \label{sec:SplitPlucker1} it follows that when $A_1>0$, $A_2<0$, and $A_1,A_2$ are odd 
\[s(\gamma) = \left(\frac{A_1/D}{-A_2/D}\right)\left(\frac{B_1/D_1}{A_1/D}\right)\left(\frac{-B_2/D_2}{-A_2/D}\right)\left(\frac{D_1}{C_1}\right)\left(\frac{D_2}{C_2}\right).\]
\end{comment}

%%%%%%%%%%%%%%%%%%%%%%%%%%%%%%%
%%%%%%%%%SPLITTING: Reduced CASE
%%%%%%%%%%%%%%%%%%%%%%%%%%%%%%%

\subsection{The Splitting: The Reduced Case} \label{sec:SplitPlucker3}

This section contains the computation of the splitting when $(A_1,A_2) = 2^\ell$. By combining this formula with the reduction step, the transformation of the formula for the splitting in terms of block parameters into a formula in terms of Pl\"{u}cker coordinates will be complete.

If the reduction step could be applied with $D=(A_1,A_2)$, the computation would be more straightforward. Unfortunately, when $2$ divides $(A_1,A_2)$, it is not a valid choice in the reduction step. So the reduction step of Subsection \ref{sec:SplitPlucker2} can only be used to remove the odd part of the GCD of $A_1$ and $A_2$. So suppose $D=(A_1,A_2) = 2^\ell.$ By applying the symmetries of Subsection \ref{sec:PluckerSym} we can impose the additional constraints $A_1>0$, and $A_2/2^{\ell}\equiv 1\,(\text{mod } 2)$; the calculations of Subsection \ref{sec:SymSpl} describe the effect of these symmetries on the value of the splitting. Specifically, the Cartan involution composed with conjugation by the long element can be used to swap $A_1$ and $A_2$; thus, if $A_2$ is more even than $A_1$ they can be swapped; if the new $A_1$ is negative, then conjugation by $t(1,1,-1)$ will flip the sign of $A_1$. Thus the computation in the general case can be reduced to that of the special case just described. The evaluation of the splitting in this special case is the content of the following proposition.

 \begin{proposition}\label{splitreduced}Let $\gamma\in\Gamma_{1}(4)$ with Pl\"{u}cker coordinates $(4A_1,4B_1,C_1,4A_2,4B_2,C_2)$ such that $A_1>0$, and $A_2/(A_1,A_2)\equiv 1\,(\text{mod } 2)$. Let $D=(A_1,A_2)=2^\ell$, $D_1 = (D,B_1)$, $D_2=D/D_1$, and let $\epsilon = \left(\frac{-1}{-B_1/D_1}\right)$. Then
\begin{equation} s(\gamma)=\label{plucksplitting}
\left(\frac{-\epsilon}{-A_1A_2}\right)\left(\frac{A_1/D}{A_2/D}\right)
\left(\frac{B_1/D_1}{A_1/D}\right)\left(\frac{4B_2/D_2}{\text{sign}(A_2)A_2/D}\right)\left(\frac{D_1}{C_1}\right)\left(\frac{D_2}{C_2}\right).\end{equation}
 \end{proposition}

\textbf{Proof:} Let $\gamma\in \Gamma_1(4)$ have Pl{\"u}cker Coordinates $(4A_1,4B_1,C_1,4A_2,4B_2,C_2)$.  Let $B_j = 2^{\delta_j}b_j$, where $j=1,2$, and $(2,b_j) = 1.$  Let $\epsilon = \left(\frac{-1}{-b_1}\right)=\left(\frac{-1}{-B_1}\right)$. 

The splitting and the extra conditions on $A_1$ and $A_2$ are right $\Gamma_{\infty}$-invariant. Thus multiplication on the right by an element of $\Gamma_{\infty}$ can be used to reduce the computation to that of a representative with still more favorable properties.  If $x,z\in\mathbb{Z}$, then $\gamma n(2^{\delta_1+2}x,0,z)$ has Plucker Coordinates $(4A_1,4\mathcal{B}_1,\mathcal{C}_1,4A_2,4\mathcal{B}_2,\mathcal{C}_2)$ where 
\begin{align}\mathcal{B}_1 =& B_1-(2^{\delta_1+2}x)A_1,\label{B1}\\
\mathcal{B}_2 =& B_2,\text{ and}\label{B2}\\
\mathcal{C}_2 =& C_2-4B_2(2^{\delta_1+2}x)-4A_2z.\label{C2}\end{align}
The factor of $2^{\delta_1+2}$ in front of $x$ ensures that ord$_2(\mathcal{B}_1)=\text{ord}_2(B_1)$ and $\mathcal{B}_1\equiv B_1\text{(mod }4)$. Now we will specify certain congruences conditions  that will aid in our computation of the splitting.
Using the Chinese Remainder Theorem, choose $x\in\mathbb{Z}$ subject to the conditions
\[\begin{array}{lll}
\mathcal{B}_1&\equiv&B_1 \hspace{.25cm} (\text{mod } 2^{\delta_1+2}A_1),\\
\mathcal{B}_1&\equiv&1 \hspace{.25cm} (\text{mod } A_2/2^{\ell}),\\
\mathcal{B}_1&<&0.\end{array}\]
Choose $z\in\mathbb{Z}$ subject to the conditions
\[\begin{array}{lll}\mathcal{C}_2&\equiv&1 \hspace{.25cm} (\text{mod } \mathcal{B}_1/2^{\delta_1}),\\
\mathcal{C}_2&\equiv&C_2-4B_22^{\delta_1+2}x \hspace{.25cm} (\text{mod } 4A_2),\\
\mathcal{C}_2&<&0.\end{array}\]
This more suitable representative will be used to compute the splitting.
%%%%%NOTE: The choice of x makes \mathcal{B}_1/2^{\delta_1} and $A_2$ relatively prime.

By Proposition \ref{blockrep}, $\gamma n(2^{\delta_1+2}x,0,z)$ can be factored so that it has block parameters
\begin{equation}\label{redblockparams}\begin{array}{lllllllll}
c_2&=&-4A_1&,&d_2&=&\pm(\mathcal{B}_1,\mathcal{C}_1)&=&\pm (\mathcal{B}_1, \frac{A_1\mathcal{C}_2+4\mathcal{B}_1\mathcal{B}_2}{-A_2}),\\
c_3&=&-4\mathcal{B}_1/d_2&,&d_3&=&-\mathcal{C}_1/d_2&=&\frac{A_1\mathcal{C}_2+4\mathcal{B}_1\mathcal{B}_2}{A_2d_2},\\
d_1&=&-\mathcal{C}_2/d_2&,&\frac{-\mathcal{B}_1}{d_2}c_1&=&-A_2+\frac{-\mathcal{C}_2}{d_2}(-A_1)(a_3),
\end{array}\end{equation}
where $a_3>0$.  Note that $(\mathcal{B}_1,A_2/2^{\ell})=(\mathcal{B}_1,\mathcal{C}_2)=1$.  Therefore $d_2=1$ as $(A_1,\mathcal{B}_1,\mathcal{C}_1)=1$.

Since $c_2,c_3\neq 0$ and $a_3>0$, equation (\ref{nonarithsplitting}) shows that 
\[s_{\text{NA}}(\gamma) = \begin{cases}
(c_1(-A_1)(-A_2),c_1)= (-A_1A_2,c_1)&, c_1\neq0\\
1&, c_1=0.\end{cases}\] 
By considering the cases $A_2>0$ and $A_2<0$ it can be shown that $s_{\text{NA}}(\gamma)=-\text{sign}(A_2).$

The computation of the arithmetic part remains. The arithmetic part of the splitting is given by $\left(\frac{-A_1}{d_2}\right)\left(\frac{-\mathcal{B}_1/d_2}{-\mathcal{C}_1/d_2}\right)\left(\frac{c_1}{-\mathcal{C}_2/d_2}\right).$ Since $(\mathcal{B}_1/2^{\delta_1},A_2)=1$ and $d_2=1$ we may conclude that
\begin{multline}
\left(\frac{-A_1}{d_2}\right)\left(\frac{-\mathcal{B}_1/d_2}{-\mathcal{C}_1/d_2}\right)\left(\frac{c_1}{-\mathcal{C}_2/d_2}\right)\\
=\left(\frac{\epsilon2^{\delta_1}}{-\mathcal{C}_1}\right)\left(\frac{-\epsilon\mathcal{B}_1/2^{\delta_1}}{\text{sign}(A_2)[(A_1/2^{\ell})\mathcal{C}_2+4(\mathcal{B}_1\mathcal{B}_2/2^{\ell})]}\right)\left(\frac{-\epsilon\mathcal{B}_1/2^{\delta_1}}{\text{sign}(A_2)A_2/2^{\ell}}\right)\left(\frac{c_1}{-\mathcal{C}_2}\right).\label{spliteqn1}\end{multline}
Now we apply the fourth property in Proposition \ref{Kronecker} and use the equation involving $c_1$ in line (\ref{redblockparams}) to get
\begin{multline}(\ref{spliteqn1})=\left(\frac{\epsilon2^{\delta_1}}{-\mathcal{C}_1}\right)\left(\frac{-\epsilon\mathcal{B}_1/2^{\delta_1}}{\text{sign}(A_2)(A_1/2^{\ell})\mathcal{C}_2}\right)\left(\frac{-\epsilon\mathcal{B}_1/2^{\delta_1}}{\text{sign}(A_2)A_2/2^{\ell}}\right)\left(\frac{\mathcal{B}_1A_2-\mathcal{C}_2A_1a_3}{-\mathcal{C}_2}\right).\label{spliteqn2}
\end{multline}
The next equality follows as $-\mathcal{C}_1\equiv 1\,(\text{mod }4)$, $-\mathcal{C}_2>0$, $\mathcal{B}_1<0$, and from an application of quadratic reciprocity. 
\begin{flalign}(\ref{spliteqn2})=\left(\frac{2^{\delta_1}}{\mathcal{C}_1}\right)\left(\frac{-\epsilon\mathcal{B}_1/2^{\delta_1}}{\text{sign}(A_2)(A_1/2^{\ell})\mathcal{C}_2}\right)\left(\frac{-\epsilon\mathcal{B}_1/2^{\delta_1}}{\text{sign}(A_2)A_2/2^{\ell}}\right)\left(\frac{-2^{\ell}\mathcal{B}_1}{\mathcal{C}_2}\right)
\left(\frac{-\mathcal{C}_2}{\text{sign}(A_2)A_2/2^{\ell}}\right).\label{spliteqn4}
\end{flalign}
Now replace $\mathcal{C}_2$ using the equation $A_1\mathcal{C}_2+4\mathcal{B}_1\mathcal{B}_2+\mathcal{C}_1A_2=0$ and the third property in Proposition \ref{Kronecker} to see that
\begin{multline}(\ref{spliteqn4})=\left(\frac{2^{\delta_1}}{\mathcal{C}_1}\right)\left(\frac{-\epsilon\mathcal{B}_1/2^{\delta_1}}{\text{sign}(A_2)(A_1/2^{\ell})\mathcal{C}_2}\right)\left(\frac{-\epsilon\mathcal{B}_1/2^{\delta_1}}{\text{sign}(A_2)A_2/2^{\ell}}\right)\left(\frac{-2^{\ell}\mathcal{B}_1}{\mathcal{C}_2}\right)\\
\times\left(\frac{A_1/2^{\ell}}{\text{sign}(A_2)A_2/2^{\ell}}\right)\left(\frac{4\mathcal{B}_1\mathcal{B}_2/2^{\ell}}{\text{sign}(A_2)A_2/2^{\ell}}\right).\label{spliteqn5}\end{multline}
Next rearrange the terms, use that $\mathcal{B}_1<0$ and $A_1>0$, and  group the terms involving $\mathcal{C}_2$ to get
\begin{align}
(\ref{spliteqn5})=&\left(\frac{A_1/2^{\ell}}{A_2/2^{\ell}}\right)\left(\frac{-\epsilon\mathcal{B}_1/2^{\delta_1}}{\text{sign}(A_2)(A_1/2^{\ell})}\right)\left(\frac{(4\mathcal{B}_1\mathcal{B}_2/2^{\ell})(-\epsilon\mathcal{B}_1/2^{\delta_1})}{\text{sign}(A_2)A_2/2^{\ell}}\right)\left(\frac{2^{\delta_1}}{\mathcal{C}_1}\right)\left(\frac{\epsilon2^{l+\delta_1}}{\mathcal{C}_2}\right).\label{place2}
\end{align}
Using basic properties of the Kronecker Symbol it follows that
\begin{align}
(\ref{place2})=&\left(\frac{\epsilon}{\text{sign}(A_2)\mathcal{C}_2}\right)\left(\frac{A_1/2^{\ell}}{A_2/2^{\ell}}\right)\left(\frac{-\epsilon\mathcal{B}_1/2^{\delta_1}}{A_1/2^{\ell}}\right)\left(\frac{(-\epsilon\mathcal{B}_2/2^{\delta_2})(2^{\delta_1+\delta_2+2-l})}{\text{sign}(A_2)A_2/2^{\ell}}\right)\left(\frac{2^{\delta_1}}{\mathcal{C}_1}\right)\left(\frac{2^{l+\delta_1}}{\mathcal{C}_2}\right).\label{spliteqn7}
\end{align}
Now we show that the formula remains the same when the modified Pl\"{u}cker coordinates are switched back to the original Pl\"{u}cker coordinates using the equations of lines (\ref{B1}), (\ref{B2}), and (\ref{C2}). The switch for $B_1$ and $B_2$ is direct, but it appears that some care is needed when dealing with $C_1$ and $C_2$. The case when $\ell>0$ is straightforward as $\mathcal{C}_i\equiv C_i\text{(mod }8)$ for $i=1,2$. Unfortunately, this need not be true when $\ell=0$. 

Supposing that $\ell=0$ there will be two cases to consider: $\delta_1=0$ and $\delta_1\neq0$.  If $\delta_1=0$ both exponents will be 0 and there is nothing to show. The last case to consider is when $\ell=0$ and $\delta_1>0$. In this case $A_1C_2+C_1A_2\equiv 0 (\text{mod }8)$ and $A_1\mathcal{C}_2+\mathcal{C}_1A_2\equiv 0 (\text{mod }8)$. As $A_1,A_2$ are odd, $C_1C_2\equiv \mathcal{C}_1\mathcal{C}_2\text{(mod }8)$. Thus $\left(\frac{2^{\delta_1}}{C_1}\right)\left(\frac{2^{\delta_1}}{C_2}\right)=\left(\frac{2^{\delta_1}}{\mathcal{C}_1}\right)\left(\frac{2^{\delta_1}}{\mathcal{C}_2}\right)$. Therefore,
\begin{multline}(\ref{spliteqn7})=\left(\frac{\epsilon}{\text{sign}(A_2)C_2}\right)\left(\frac{A_1/2^{\ell}}{A_2/2^{\ell}}\right)\left(\frac{-\epsilon B_1/2^{\delta_1}}{A_1/2^{\ell}}\right)\left(\frac{(-\epsilon B_2/2^{\delta_2})2^{\delta_1+\delta_2+2-l}}{\text{sign}(A_2)A_2/2^{\ell}}\right)\left(\frac{2^{\delta_1}}{C_1}\right)\left(\frac{2^{l+\delta_1}}{C_2}\right)\label{quadprod}.\end{multline}

Again the computation will bifurcate as we try to simplify the last four Kronecker Symbols that appear in the previous line.  Case 1 will consist of $\delta_1\leq l$ and Case 2 will consist of $\delta_1>l.$ Case 1 is the easier of the two, so we will focus on Case 2.
\begin{comment}
Case 1:
\noindent If $\delta_1\leq l$ then $D_1=(2^{\ell},B_1) = 2^{\delta_1}$ and $D_2 = 2^{l-\delta_1}$.  Thus,
\begin{align*}\left(\frac{-B_1/2^{\delta_1}}{A_1/2^{\ell}}\right) &= \left(\frac{-B_1/D_1}{A_1/2^{\ell}}\right),\\
\left(\frac{-\epsilon B_2/2^{\delta_2}}{\text{sign}(A_2)A_2/2^{\ell}}\right)\left(\frac{2^{\delta_1+\delta_2+2-l}}{\text{sign}(A_2)A_2/2^{\ell}}\right)&=\left(\frac{-\epsilon 4B_2/D_2}{\text{sign}(A_2)A_2/2^{\ell}}\right),\\
\left(\frac{2^{\delta_1}}{C_1}\right)\left(\frac{2^{l+\delta_1}}{C_2}\right)=\left(\frac{2^{\delta_1}}{C_1}\right)\left(\frac{2^{l-\delta_1}}{C_2}\right)&=\left(\frac{D_1}{C_1}\right)\left(\frac{D_2}{C_2}\right).\end{align*}
\end{comment}

\noindent Case 2:
\noindent If $\delta_1> l$ then $D_1=(2^{\ell},B_1) = 2^{\ell}$ and $D_2 = 1$.  Thus,
\begin{multline*}
\left(\frac{-B_1/2^{\delta_1}}{A_1/2^{\ell}}\right)\left(\frac{-\epsilon B_2/2^{\delta_2}}{\text{sign}(A_2)A_2/2^{\ell}}\right)\left(\frac{2^{\delta_1+\delta_2+2-l}}{\text{sign}(A_2)A_2/2^{\ell}}\right)\left(\frac{2^{\delta_1}}{C_1}\right)\left(\frac{2^{l+\delta_1}}{C_2}\right) \\= \left(\frac{-B_1/D_1}{A_1/2^{\ell}}\right)\left(\frac{-\epsilon 4B_2/D_2}{\text{sign}(A_2)A_2/2^{\ell}}\right)\left(\frac{D_1}{C_1}\right)\left(\frac{D_2}{C_2}\right)\left(\frac{2^{\delta_1-l}}{C_2A_1/2^{\ell}}\right)\left(\frac{2^{\delta_1-l}}{C_1(-A_2)/2^{\ell}}\right).\end{multline*}
Finally we must show that $\left(\frac{2^{\delta_1-l}}{C_2A_1/2^{\ell}}\right)\left(\frac{2^{\delta_1-l}}{C_1(-A_2)/2^{\ell}}\right) = 1$.
 
To see this note that $\delta_1>l$ implies that $\delta_1>0$. Thus $(A_1/2^{\ell})C_2+C_1(A_2/2^{\ell})\equiv 0 \text{ }(\text{mod } 8)$ or equivalently, $(A_1/2^{\ell})C_2\equiv-C_1(A_2/2^{\ell}) \text{ }(\text{mod } 8)$. Thus $\left(\frac{2^{\delta_1-l}}{C_2A_1/2^{\ell}}\right)\left(\frac{2^{\delta_1-l}}{C_1(-A_2)/2^{\ell}}\right) = 1$.

In both Case 1 and Case 2 the arithmetic part of the splitting is given by 
\[s_A(\gamma) = \left(\frac{\epsilon}{\text{sign}(A_2)C_2}\right)\left(\frac{A_1/D}{A_2/D}\right)\left(\frac{-\epsilon B_1/D_1}{A_1/D}\right)\left(\frac{-\epsilon 4B_2/D_2}{\text{sign}(A_2)A_2/D}\right)\left(\frac{D_1}{C_1}\right)\left(\frac{D_2}{C_2}\right).\]
After incorporating the non-arithmetic factor and groupings the $\epsilon$ terms and the $-1$ terms together the formula simplifies to 
\[s(\gamma) = \left(\frac{-\epsilon}{-A_1A_2}\right)\left(\frac{A_1/D}{A_2/D}\right)\left(\frac{B_1/D_1}{A_1/D}\right)\left(\frac{4B_2/D_2}{\text{sign}(A_2)A_2/D}\right)\left(\frac{D_1}{C_1}\right)\left(\frac{D_2}{C_2}\right).\]\EndProof

At this point we can combine Proposition \ref{splitreduction} and Proposition \ref{splitreduced} to prove Theorem \ref{theorem:PluckSplit}.

\textbf{Proof of Theorem \ref{theorem:PluckSplit}:} Let $\gamma\in\Gamma_{1}(4)$ with Pl\"{u}cker coordinates \newline$(4A_1,4B_1,C_1,4A_2,4B_2,C_2)$ such that $A_1>0$, and $A_2/(A_1,A_2)\equiv 1\,(\text{mod } 2)$. Let $D=(A_1,A_2)$, $D_1 = (D,B_1)$, $D_2=D/D_1$, and let $\epsilon = \left(\frac{-1}{-B_1/D_1}\right)$. Let $D^\prime$, $D_1^\prime$, and $D_2^\prime$ denote the odd part of $D$, $D_1$, and $D_2$ respectively. Let $S=t(1,(D_2^\prime)^{-1},(D^\prime)^{-1})$. Since $S\gamma S^{-1}$ has Pl\"{u}cker coordinates $(4A_1/D^\prime,4B_1/(D_1)^\prime,C_1,4A_2/D^\prime,4B_2/(D_2)^\prime,C_2)$, Proposition \ref{splitreduction} shows that $s(\gamma)=s(S\gamma S^{-1})\left(\frac{D_1^\prime}{C_1}\right)\left(\frac{D_2^\prime}{C_2}\right)$. Let $\epsilon^\prime=\left(\frac{-1}{-B_1/(D_1^\prime)}\right)$ and note that $\epsilon^\prime=\epsilon$. Since $(A_1/D^\prime,A_2/D^\prime)$ is a power of $2$ we can apply Proposition \ref{splitreduced} to evaluate $s(S\gamma S^{-1})$. In particular, we have 
\begin{align*}s(S\gamma S^{-1}) =&\left(\frac{-\epsilon^\prime}{-A_1A_2}\right)\left(\frac{A_1/D}{A_2/D}\right)\left(\frac{B_1/D_1}{A_1/D}\right)\left(\frac{4B_2/D_2}{\text{sign}(A_2)A_2/D}\right)\left(\frac{(D_1/D_1^\prime)}{C_1}\right)\left(\frac{(D_2/D_2^\prime)}{C_2}\right).
\end{align*} 
As $s(\gamma)=s(S\gamma S^{-1})\left(\frac{D_1^\prime}{C_1}\right)\left(\frac{D_2^\prime}{C_2}\right)$ and $\epsilon=\epsilon^\prime$ we can conclude that
\[s(\gamma) = \left(\frac{-\epsilon}{-A_1A_2}\right)\left(\frac{A_1/D}{A_2/D}\right)\left(\frac{B_1/D_1}{A_1/D}\right)\left(\frac{4B_2/D_2}{\text{sign}(A_2)A_2/D}\right)\left(\frac{D_1}{C_1}\right)\left(\frac{D_2}{C_2}\right).\]\EndProof

%%%%%%%%%%%%%%%%%%%%%%%%%%%%%%%%%%%%
%%%Twisted Multiplicativity
%%%%%%%%%%%%%%%%%%%%%%%%%%%%%%%%%%%%

\subsection{Twisted Multiplicativity} \label{sec:Twist}

\begin{proposition} \label{twistmult} Let $A_1,\alpha_1\in\mathbb{Z}_{>0}$, $A_2,\alpha_2\in \mathbb{Z}$ such that $A_1,A_2$ are odd, \newline$(A_1A_2,\alpha_1\alpha_2)=1$, $A_1\alpha_1+A_2\alpha_2\equiv 0\text{ (mod }4)$, and $\frac{\alpha_2}{(\alpha_1,\alpha_2)}\equiv 1\text{ (mod }2)$. Let $\mu = \left(\frac{-1}{-A_1A_2}\right)$. Then with respect to the map from Proposition \ref{mult}, \newline$\phi:S(A_1\alpha_1,A_2\alpha_2)\rightarrow S(A_1,\mu A_2)\times S(\alpha_1,-\mu\alpha_2)$, the following holds:
\begin{equation}s(\gamma) = s(\pi_1(\phi(\gamma)))s(\pi_2(\phi(\gamma)))\left(\frac{\alpha_2}{(\frac{-1}{A_1})A_1}\right)\left(\frac{\alpha_1}{A_2}\right),\end{equation}
where $\pi_i$ is projection onto the $i$-th factor.
\end{proposition}

Note that the visible asymmetry in the formula is a result of the asymmetry contained in the hypotheses. In particular, since $\alpha_1>0$, $\left(\frac{\alpha_1}{\pm A_2}\right)=\left(\frac{\alpha_1}{A_2}\right).$

\noindent\textbf{Proof:} Recall the map
 \begin{multline}(A_1\alpha_1,B_1,C_1,A_2\alpha_2,B_2,C_2)
 \stackrel{\phi}{\mapsto}((A_1,B_1,\mathcal{C}_1,\mu A_2,B_2,\gamma C_2),\\
 (\alpha_1,B_1,\left(\frac{-1}{A_2}\right) A_2C_1,-\mu\alpha_2,-\left(\frac{-1}{A_1}\right)B_2,\left(\frac{-1}{A_1}\right)A_1C_2)).
 \end{multline}
Let $D=(A_1\alpha_1,A_2\alpha_2)$, $d=(A_1,A_2)$, and $\delta=(\alpha_1,\alpha_2)$. Then $D=(A_1\alpha_1,A_2\alpha_2) = (A_1,A_2)(\alpha_1,\alpha_2) = d\delta$, as $(A_i,\alpha_j) =1$. Similarly, define $D_1 = (D,B_1)$,  $d_1=(d,B_1)$, and $\delta_1 = (\delta,B_1)$. Again $D_1 = d_1\delta_1$. Let $D_2 = \frac{D}{D_1}$, $d_2 = \frac{d}{d_1}$, and $\delta_2 = \frac{\delta}{\delta_2}$. Note that $D_2=d_2\delta_2$. Let $\epsilon = \left(\frac{-1}{-B_1/D_1}\right)$, $\epsilon_1 = \left(\frac{-1}{-B_1/d_1}\right)$, and $\epsilon_2=\left(\frac{-1}{-B_1/\delta_1}\right)$.

Begin with the formula of the splitting described in Theorem \ref{theorem:PluckSplit}, since $A_1+\mu A_2\equiv 0\, (\text{mod } 4)$ it follows that
\begin{equation*}
s(\gamma)=\left(\frac{-\epsilon}{\alpha_1\mu\alpha_2}\right)\left(\frac{A_1\alpha_1/D}{A_2\alpha_2/D}\right)\left(\frac{B_1/D_1}{A_1\alpha_1/D}\right)\left(\frac{4B_2/D_2}{\text{sign}(A_2\alpha_2)A_2\alpha_2/D}\right)\left(\frac{D_1}{C_1}\right)\left(\frac{D_2}{C_2}\right).\end{equation*}
Similarly,
\begin{align*}
s(\pi_1(\phi(\gamma)))=&\left(\frac{A_1/d}{\mu A_2/d}\right)\left(\frac{B_1/d_1}{A_1/d}\right)\left(\frac{4B_2/d_2}{\text{sign}(\mu A_2)\mu A_2/d}\right)\left(\frac{d_1}{\mathcal{C}_1}\right)\left(\frac{d_2}{\gamma C_2}\right),\end{align*}
and
\[s(\pi_2(\phi(\gamma)))=\left(\frac{-\epsilon_2}{\mu\alpha_1\alpha_2}\right)\left(\frac{\alpha_1/\delta}{-\mu\alpha_2/\delta}\right)\left(\frac{B_1/\delta_1}{\alpha_1/\delta}\right)\left(\frac{-4\left(\frac{-1}{A_2}\right)\mu B_2/\delta_2}{\text{sign}(\alpha_2)\alpha_2/\delta}\right)\left(\frac{\delta_1}{ A_2C_1}\right)\left(\frac{\delta_2}{A_1C_2}\right)\]

Next we consider the product of the analogous terms in the formula of each splitting. First consider the terms involving $\epsilon$.
\begin{align}
\left(\frac{-\epsilon}{\alpha_1\mu\alpha_2}\right)\left(\frac{-\epsilon_2}{\mu\alpha_1\alpha_2}\right)=\left(\frac{\epsilon\epsilon_2}{\alpha_1\mu\alpha_2}\right)
=\left(\frac{\left(\frac{-1}{d_1}\right)}{\alpha_1\mu\alpha_2}\right).\label{twist1}\end{align}
Consider the terms involving $A_1$, $A_2$, $\alpha_1$, or $\alpha_2$. We can simplify this expression to get
\begin{equation}
\left(\frac{A_1\alpha_1/D}{A_2\alpha_2/D}\right)\left(\frac{A_1/d}{\mu A_2/d}\right)\left(\frac{\alpha_1/\delta}{-\mu\alpha_2/\delta}\right)
\label{twist2}=\left(\frac{A_1/d}{(\alpha_2/\delta)}\right)\left(\frac{\alpha_1/\delta}{A_2/d}\right).
\end{equation}
Next we simplify the terms involving $B_1$ to get
\begin{equation}
\left(\frac{B_1/D_1}{A_1\alpha_1/D}\right)\left(\frac{B_1/d_1}{A_1/d}\right)\left(\frac{B_1/\delta_1}{\alpha_1/\delta}\right)
\label{twist3}=\left(\frac{\delta_1}{A_1/d}\right)\left(\frac{d_1}{\alpha_1/\delta}\right).
\end{equation}
Now we consider the terms involving $B_2$.
\begin{multline}
\left(\frac{4B_2/D_2}{\text{sign}(A_2\alpha_2)A_2\alpha_2/D}\right)\left(\frac{4B_2/d_2}{\text{sign}(\mu A_2)\mu A_2/d}\right)\left(\frac{-4\left(\frac{-1}{A_2}\right)\mu B_2/\delta_2}{\text{sign}(\alpha_2)\alpha_2/\delta}\right)\\
\label{twist4}=\left(\frac{\delta_2}{A_2/d}\right)\left(\frac{-\mu\left(\frac{-1}{A_2}\right)d_2}{\text{sign}(\alpha_2)\alpha_2/\delta}\right).
\end{multline}
Finally, using quadratic reciprocity we can simplify the terms involving $C_1$, $\mathcal{C}_1$, or $C_2$. 
\begin{equation}
\left(\frac{D_1}{C_1}\right)\left(\frac{d_1}{\mathcal{C}_1}\right)\left(\frac{\delta_1}{A_2C_1}\right)\left(\frac{D_2}{C_2}\right)\left(\frac{d_2}{\gamma C_2}\right)\left(\frac{\delta_2}{A_1C_2}\right)
\label{twist5}=\left(\frac{\mu\alpha_2}{d_1}\right)\left(\frac{\delta_1}{A_2}\right)\left(\frac{\alpha_1}{d_2}\right)\left(\frac{\delta_2}{A_1}\right).\end{equation}

Now the task is to simplify the product of lines (\ref{twist1})-(\ref{twist5}). After rearranging the terms the quantity is given by 
\begin{multline}\label{twistnew1}\bigg(\frac{(-1)^{(d_1-1)/2}}{\alpha_1\mu\alpha_2}\bigg)\bigg(\frac{A_1/d}{\alpha_2/\delta}\bigg)\bigg(\frac{\alpha_1/\delta}{A_2/d}\bigg)\bigg(\frac{\delta_2}{A_2/d}\bigg)\bigg(\frac{\delta_1}{A_1/d}\bigg)\bigg(\frac{\delta_1}{A_2}\bigg)\\
\times\bigg(\frac{d_1}{\alpha_1/\delta}\bigg)\bigg(\frac{\mu\alpha_2}{d_1}\bigg)\bigg(\frac{\delta_2}{A_1}\bigg)\bigg(\frac{(-1)^{(A_1-1)/2}d_2}{\text{sign}(\alpha_2)\alpha_2/\delta}\bigg)\bigg(\frac{\alpha_1}{d_2}\bigg).
\end{multline}

A lengthy but elementary calculation shows that 

\begin{align}(\ref{twistnew1}) \nonumber=&\bigg(\frac{(-1)^{(d_1-1)/2}}{\alpha_1\mu\alpha_2}\bigg)\bigg(\frac{A_1/d_1}{\alpha_2/\delta}\bigg)((-1)^{(A_1-1)/2},\text{sign}(\alpha_2))\bigg(\frac{\alpha_1}{A_2/d_1}\bigg)\\
&\hspace{3cm}\times\bigg(\frac{\delta}{A_1}\bigg)\bigg(\frac{(-1)^{(\alpha_1/\delta-1)/2}\alpha_1/\delta}{d_1}\bigg)\bigg(\frac{(-1)^{(A_1-1)/2}}{\alpha_2/\delta}\bigg)\bigg(\frac{\mu\alpha_2}{d_1}\bigg)\\
\nonumber=&\bigg(\frac{(-1)^{(d_1-1)/2}}{\alpha_1\mu\alpha_2}\bigg)\bigg(\frac{\alpha_2}{(-1)^{(A_1-1)/2}A_1}\bigg)\bigg(\frac{\alpha_1}{A_2}\bigg)\bigg(\frac{(-1)^{(\alpha_1/\delta_1-1)/2}}{d_1}\bigg)
\bigg(\frac{(-1)^{(d_1-1)/2}}{\alpha_2/\delta}\bigg)\bigg(\frac{\mu}{d_1}\bigg)\\
=&\bigg(\frac{\alpha_2}{(-1)^{(A_1-1)/2}A_1}\bigg)\bigg(\frac{\alpha_1}{A_2}\bigg).
\end{align}

\EndProof

%%%%%%%%%%%%%%%%%%%%%%%%%%%%%%%%%%%%%%%%%%%%%%%%%%%%%
%%%%%%%%%THE SPLITTING ON OTHER CELLS
%%%%%%%%%%%%%%%%%%%%%%%%%%%%%%%%%%%%%%%%%%%%%%%%%%%%%

\subsection{The Splitting on Other Cells} \label{sec:SplitNotBig}

The content of Subsection \ref{sec:PluckSplit} is a description of the splitting map on the big cell; this section collects the description of $s$ on the smaller Bruhat cells.

\begin{proposition} Let $\gamma\in \Gamma_{1}(4)$ with Pl\"{u}cker Coordinates $(A_1,B_1,C_1,A_2,B_2,C_2)$. Then:

\begin{center}\begin{tabular}{|c|c|c|} \hline
\text{Cell} & $(A_1,B_1,C_1,A_2,B_2,C_2)$ & $s(\gamma)$\\ \hline
$B$ & $(0,0,-1,0,0,-1)$ & 1\\ \hline
$Bw_{\alpha_{1}}B$ & $(0,0,-1,0,B_2,C_2)$ & $\left(\frac{B_2}{-C_2}\right)$ \\ \hline
$Bw_{\alpha_{2}}B$ & $(0,B_1,C_1,0,0,-1)$ & $\left(\frac{-B_1}{-C_1}\right)$\\ \hline
$Bw_{\alpha_{1}}w_{\alpha_{2}}B$ & $(0,B_1,C_1,A_2,B_2,C_2)$ & $\left(\frac{A_2/B_1}{-C_2}\right)\left(\frac{-B_1}{-C_1}\right)$\\ \hline
$Bw_{\alpha_{2}}w_{\alpha_{1}}B$ & $(A_1,B_1,C_1,0,B_2,C_2)$ & $(-A_1,B_2)\left(\frac{-A_1/B_2}{-C_1}\right)\left(\frac{B_2}{-C_2}\right)$\\ \hline
$Bw_{\ell}B$ & $(A_1,B_1,C_1,A_2,B_2,C_2)$ & Equation (\ref{plucksplitting}) \\ \hline
\end{tabular}
\end{center}
\end{proposition}

\noindent\textbf{Proof:} We will focus on just two of the cases; the remaining cases are similar.

\noindent ($Bw_{\alpha_1}B$) In this case Proposition \ref{blockPluck} implies that $c_2=c_3=0$, $d_2=d_3=1$, $c_1=B_2$, and $d_1=-C_2$.
Thus $\gamma=n\left(\begin{smallmatrix}
a_1 & b_1 & 0\\
B_2 & -C_2 & 0\\
0 & 0 & 1
\end{smallmatrix}\right)n^\prime$,
where $n,n^\prime\in \Gamma_{\infty}.$ As $c_2=0$ the nonarithmetic factor will be equal to 1, by equation (\ref{nonarithsplitting}). Thus by the left and right $\Gamma_{\infty}$-invariance of $s$ and equation (\ref{blocksplitting}), $s(\gamma)=\left(\frac{B_2}{-C_2}\right)$.

\noindent ($Bw_{\alpha_2}w_{\alpha_1}B$) This case will utilize the symmetry of $s$ with respect to the Cartan involution composed with the conjugation by the long element.  In section \ref{sec:SymSpl} we showed that $s(\phi(\gamma)) = (-A_1,B_2)s(\gamma).$  By combining this identity with the result for $Bw_{\alpha_{1}}w_{\alpha_{2}}B$ it follows that $s(\gamma) = (-A_1,B_2)\left(\frac{-A_1/B_2}{-C_1}\right)\left(\frac{B_2}{-C_2}\right)$.

%\bibliographystyle{plain}
%\bibliography{ThesisBib}

\begin{thebibliography}{9}

\bibitem{BLS99} 
Banks, William D.; Levy, Jason; Sepanski, Mark R.
Block-Compatible Metaplectic Cocycles,
\textit{Journal f{\"u}r die reine und angewandte Mathematik},
Volume: 507 (1999) 

\bibitem{BGL3} 
Bump, Daniel.
\textit{Automorphic Forms on GL$(3,\mathbb{R})$}. 
Springer-Verlag, New York, 1984.

\bibitem{BBFH07}
Brubaker, B.; Bump, D.; Friedberg, S.; Hoffstein, J.,
Weyl group multiple Dirichlet series. III. Eisenstein series and twisted unstable $A_r$.
\textit{Ann. of Math.} (2) 166 (2007), no. 1, 293-316.

\bibitem{IK04}
Iwaniec, H; Kowalski, E.,
\textit{Analytic Number Theory}.
American Mathematical Society, Rhode Island, 2004.

\bibitem{M06} 
Miller, S. D.,
Guide for the Metaplexed.
Unpublished Notes.



\end{thebibliography}

\end{document}